\documentclass[a4paper,11pt]{article}

\usepackage[latin1]{inputenc}

\usepackage{latexsym}

\usepackage{amssymb}
\usepackage{amsfonts}
\usepackage{amsmath}
\usepackage{amsthm}
\usepackage{textcomp}

\newtheorem{theorem}{Theorem}[section]    
\newtheorem{lemma}[theorem]{Lemma}          
\newtheorem{corollary}[theorem]{Corollary}
\newtheorem{proposition}[theorem]{Proposition}

\theoremstyle{definition}

\title{Logarithmic limit sets of real semi-algebraic sets}
\author{Daniele Alessandrini \\ \small \ \\ \small \textit{Universit\`a di Pisa, Italy} \\ \small \textit{E-mail address:} daniele.alessandrini@gmail.com}
\date{}

\newcommand{\nuovo}[1]{{{\bfseries \upshape #1}}}

\newcommand{\enne}{{\mathbb{N}}}
\newcommand{\ze}{{\mathbb{Z}}}
\newcommand{\qu}{{\mathbb{Q}}}
\newcommand{\erre}{{\mathbb{R}}}
\newcommand{\ci}{{\mathbb{C}}}

\newcommand{\pro}{{\mathbb{P}}}
\newcommand{\cappa}{{\mathbb{K}}}
\newcommand{\effe}{{\mathbb{F}}}

\newcommand{\ocors}{{\mathcal{O}}}

\newcommand{\ameba}{{\mathcal{A}}}

\DeclareMathOperator{\Log}{Log}
\DeclareMathOperator{\Exp}{Exp}

\DeclareMathOperator{\Span}{Span}

\newenvironment{gmatrice}{ \begin{pmatrix} }{ \end{pmatrix} }
\newenvironment{pmatrice}{ \left( \begin{smallmatrix} }{ \end{smallmatrix}  \right) }

\newcommand{\freccia}{{\ \longrightarrow \ }}

\newcommand{\tende}{{\ \rightarrow \ }}

\usepackage[dvips]{graphicx}

\newcommand{\figuresize}{4.9cm} 

\begin{document}

\sloppy

\maketitle

\begin{abstract}
This paper is about the logarithmic limit sets of real semi-algebraic sets, and, more generally, about the logarithmic limit sets of sets definable in an o-minimal, polynomially bounded structure. We prove that most of the properties of the logarithmic limit sets of complex algebraic sets hold in the real case. This include the polyhedral structure and the relation with the theory of non-archimedean fields, tropical geometry and Maslov dequantization.
\end{abstract}

\section{Introduction}

Logarithmic limit sets of complex algebraic sets have been extensively studied. They first appeared in Bergman's paper \cite{Be71}, and then they were further studied by Bieri and Groves in \cite{BG81}. Recently their relations with the theory of non-archimedean fields and tropical geometry were discovered (see for example \cite{SS04}, \cite{EKL06} and \cite{BJSST07}). They are now usually called tropical varieties, but they appeared also under the names of Bergman fans, Bergman sets, Bieri-Groves sets or non-archimedean amoebas. The logarithmic limit set of a complex algebraic set is a polyhedral complex of the same dimension as the algebraic set, it is described by tropical equations and it is the image, under the component-wise valuation map, of an algebraic set over an algebraically closed non-archimedean field. The tools used to prove these facts are mainly algebraic and combinatorial.

In this paper we extend these results to the logarithmic limit sets of real algebraic and semi-algebraic sets. The techniques we use to prove these results in the real case are very different from the ones used in the complex case. Our main tool is the cell decomposition theorem, as we prefer to look directly at the geometric set, instead of using its equations. In the real case, even if we restrict our attention to an algebraic set, it seems that the algebraic and combinatorial properties of the defining equations don't give enough information to study the logarithmic limit set.

In the following we often need to act on ${(\erre_{>0})}^n$ with maps of the form:
$$\phi_A(x_1, \dots, x_n) = (x_1^{a_{11}} \cdots x_n^{a_{1n}}, x_1^{a_{21}} \cdots x_n^{a_{2n}}, \dots, x_1^{a_{n1}} \cdots x_n^{a_{nn}} )$$
where $A = (a_{ij})$ is an $n\times n$ matrix. When the entries of $A$ are not rational, the image of a semi-algebraic set is, in general not semi-algebraic. Actually, the only thing we can say about images of semi-algebraic sets via these maps is that they are definable in the structure of the real field expanded with arbitrary power functions. This structure, usually denoted by $\overline{\erre}^\erre$, is o-minimal and polynomially bounded, and these are the main properties we need in the proofs. Moreover, if $S$ is a set definable in $\overline{\erre}^\erre$, then the image $\phi_A(S)$ is again definable, as the functions $x \freccia x^\alpha$ are definable. This property is equivalent to say that $\overline{\erre}^\erre$ has field of exponents $\erre$.

In this sense the category of semi-algebraic sets is too small for our methods. It seems that the natural context for the study of logarithmic limit sets is to fix a general expansion of the structure of the real field that is o-minimal and polynomially bounded with field of exponents $\erre$. For sets definable in such a structure, the properties that were known for the complex algebraic sets also hold. We can prove that these logarithmic limit sets are polyhedral complexes with dimension less than or equal to the dimension of the definable set, and they are the image, under the component-wise valuation map, of an extension of the definable set to a real closed non-archimedean field. An analysis of the defining equations and inequalities is carried out, showing that the logarithmic limit set of a closed semi-algebraic set can be described applying the Maslov dequantization to a suitable formula defining the semi-algebraic set. Then we show how the relation between tropical varieties and images of varieties defined over non-archimedean fields, well known for algebraically closed fields, can be extended to the case of real closed fields. We give the notion of non-archimedean amoebas of semi-algebraic sets and sets definable in other o-minimal structures and we study their relations with logarithmic limit sets of definable sets in $\erre^n$, and with patchworking families of definable sets. Note that this notion generalizes the notion of non-archimedean amoebas of semi-linear sets that have been used in \cite{DY} to study tropical polytopes.

Our motivation for this work comes from the study of Teichm\"uller spaces and, more generally, of spaces of geometric structures on manifolds. In the papers \cite{A1} and \cite{A2} we present a construction of compactification using the logarithmic limit sets. The properties of logarithmic limit sets we prove here will be used in \cite{A1} to describe the compactification. For example, the fact that logarithmic limit sets of real semi-algebraic sets are polyhedral complexes will provide an independent construction of the piecewise linear structure on the Thurston boundary of Teichm\"uller spaces. Moreover the relations with tropical geometry and the theory of non-archimedean fields will be used in \cite{A2} for constructing a geometric interpretation of the boundary points.

A brief description of the following sections. In section \ref{sez:preliminaries} we define a notion of logarithmic limit sets for general subsets of ${(\erre_{>0})}^n$, and we report some preliminary notions of model theory and o-minimal geometry that we will use in the following, most notably the notion of regular polynomially bounded structures.

In section \ref{sez:poly} we prove that logarithmic limit sets of definable sets in a regular polynomially bounded structure are polyhedral complexes with dimension less than or equal to the dimension of the definable set, and we provide a local description of these sets. The main tool we use in this section is the cell decomposition theorem.

In section \ref{sez:non-arch} we consider a special class of non-archimedean fields: the Hardy fields of regular polynomially bounded structures. These are non-archimedean real closed fields of rank one extending $\erre$, with a canonical real valued valuation and residue field $\erre$. The elements of these fields are germs of definable functions, hence they have better geometric properties than the fields of formal series usually employed in tropical geometry. The image, under the component-wise valuation map, of definable sets in the Hardy fields are related with the logarithmic limit sets of real definable sets, and with the limit of real patchworking families.

In section \ref{sez:comp} we compare the construction of this paper with other known constructions. We show that the logarithmic limit sets of complex algebraic sets are only a particular case of the logarithmic limit sets of real semi-algebraic sets, and the same happens for the limit of complex patchworking families. Hence our methods provide an alternative proof (with a topological flavor) for some known results about complex sets. We also compare the logarithmic limit sets of real algebraic sets with the construction of Positive Tropical Varieties (see \cite{SW}). Even if in many examples these two notions coincide, we show some examples where they differ.

In section \ref{sez:trop} we show how the construction of Maslov dequantization provide a relation between logarithmic limit sets of semi-algebraic sets and tropical geometry.

\section{Preliminaries}   \label{sez:preliminaries}

\subsection{Some notations}

If $x \in \erre^n$ we will denote its coordinates by $x_1, \dots, x_n$. If $\omega \in \enne^n$ we will use the multi-index notation for powers: $x^\omega = x_1^{\omega_1} \dots x_n^{\omega_n}$. We will consider also powers with real exponents, if the base is positive, hence if $x \in {(\erre_{>0})}^n$ and $\omega \in \erre^n$ we will write $x^\omega = x_1^{\omega_1} \dots x_n^{\omega_n}$.

If $s(k)$ is a sequence in $\erre^n$, we will denote its $k$-th element by $s(k) \in \erre^n$, and the coordinates by $s_i(k) \in \erre$.

Given a real number $\alpha > 1$, we will denote by $\Log_\alpha$ the \nuovo{component-wise logarithm map}, and by $\Exp_\alpha$ its inverse:
$$\Log_\alpha:{(\erre_{>0})}^n \ni (x_1, \dots, x_n) \freccia (\log_\alpha(x_1), \dots, \log_\alpha(x_n)) \in \erre^n$$
$$\Exp_\alpha: \erre^n \ni (x_1, \dots, x_n) \freccia (\alpha^{x_1}, \dots, \alpha^{x_n} ) \in {(\erre_{>0})}^n$$ 

We define a notion of \nuovo{limit} for every one-parameter family of subsets of $\erre^n$. Suppose that for all $t \in (0,\varepsilon)$ we have a set $S_t \subset \erre^n$. We can construct the deformation
$$\mathcal{D}(S_\cdot) = \{(x,t) \in \erre^n \times (0,\varepsilon) \ |\ x \in S_t\}$$
We denote by $\overline{\mathcal{D}(S)}$ the closure of $\mathcal{D}(S)$ in $\erre^n \times [0,\varepsilon)$, then we define
$$\lim_{t\tende 0} S_t = \pi( \overline{\mathcal{D}(S)} \cap \erre^n \times \{0\} ) \subset \erre^n$$
where $\pi:\erre^n \times [0,\varepsilon) \freccia \erre^n$ is the projection on the first factor.

This limit is well defined for every family of subsets of $\erre^n$.

\begin{proposition}   \label{prop:generic properties}
The set $S = \lim_{t\tende 0} S_t$ is a closed subset of $\erre^n$. A point $y$ is in $S$ if and only if there exist a sequence $y(k)$ in $\erre^n$ and a sequence $t(k)$ in $(0,\varepsilon)$ such that $t(k) \tende 0$, $y(k) \tende x$ and $\forall k \in \enne: y(k) \in S_{t(k)}$.
\end{proposition}

\subsection{Logarithmic limit sets of general sets}

Given a set $V \subset {(\erre_{>0})}^n$ and a number $t \in (0,1)$, the \nuovo{amoeba} of $V$ is
$$\ameba_t(V) = \Log_{\left(\frac{1}{t}\right)}(V) = \frac{-1}{\log_e(t)} \Log_e(V) \subset \erre^n$$
The limit of the amoebas is the \nuovo{logarithmic limit set} of $V$:
$$\ameba_0(V) = \lim_{t \tende 0} \ameba_t(V)$$

\begin{figure}[p]
\begin{center}
\includegraphics[width=\figuresize]{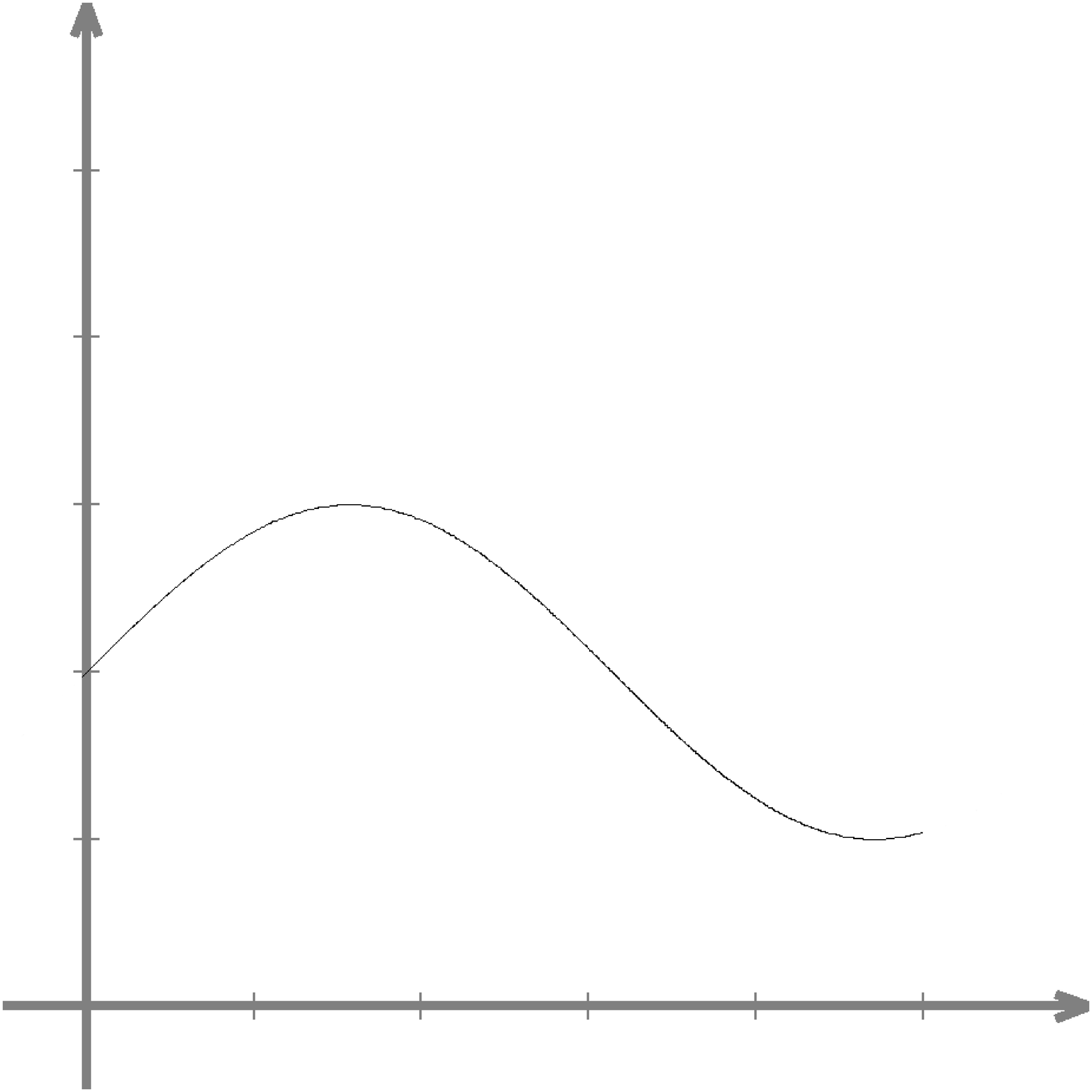}
\hspace{0.5cm}
\includegraphics[width=\figuresize]{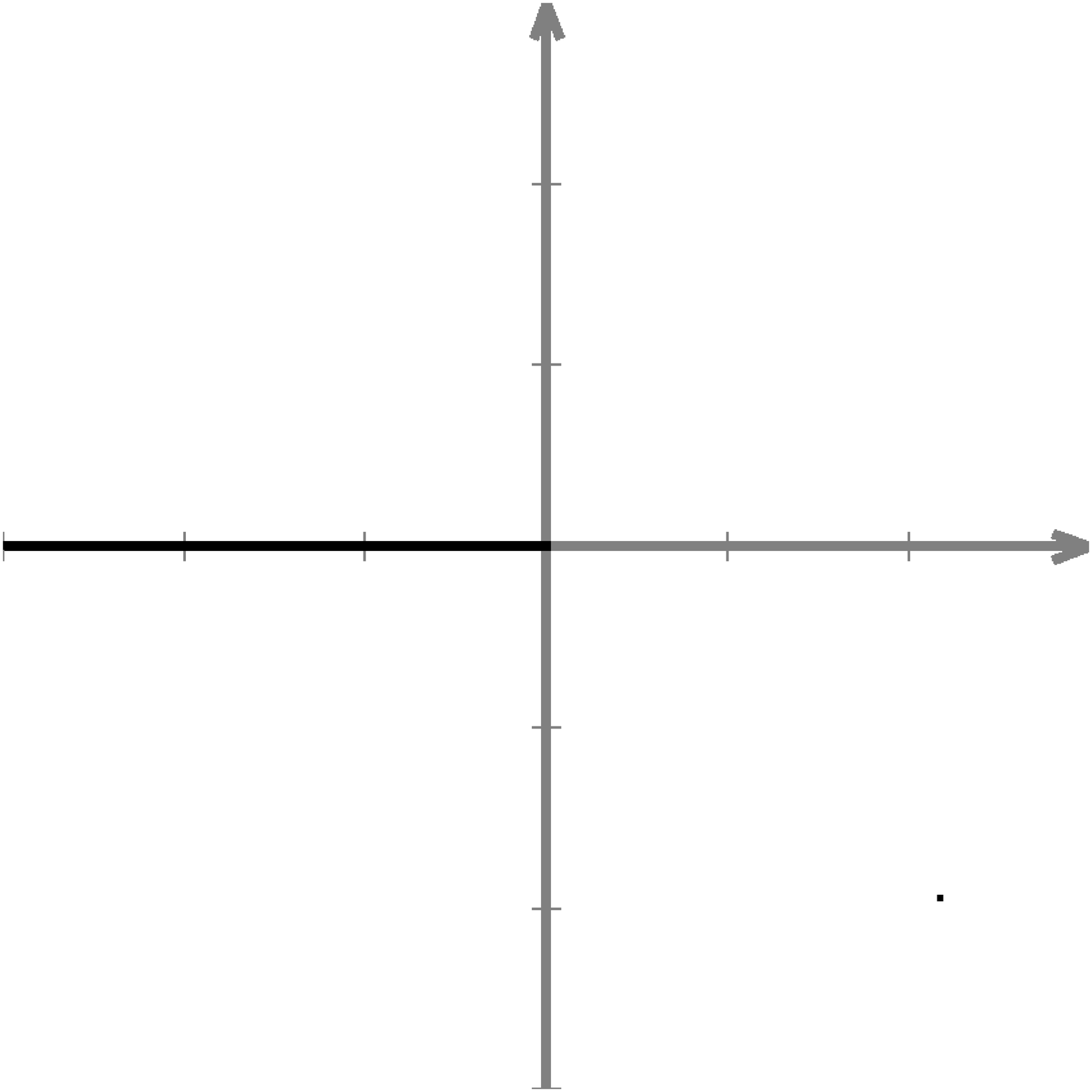}
\caption{$V = \{(x,y) \in {(\erre_{>0})}^2 \ |\ y = \sin x + 2, x \leq 5 \}$ (left picture), then $\ameba_0(V) = \{(x,y) \in \erre^2 \ |\ y = 0, x \leq 0\}$ (right picture).}  \label{fig:sin}
\end{center}
\end{figure}

\begin{figure}[p]
\begin{center}
\includegraphics[width=\figuresize]{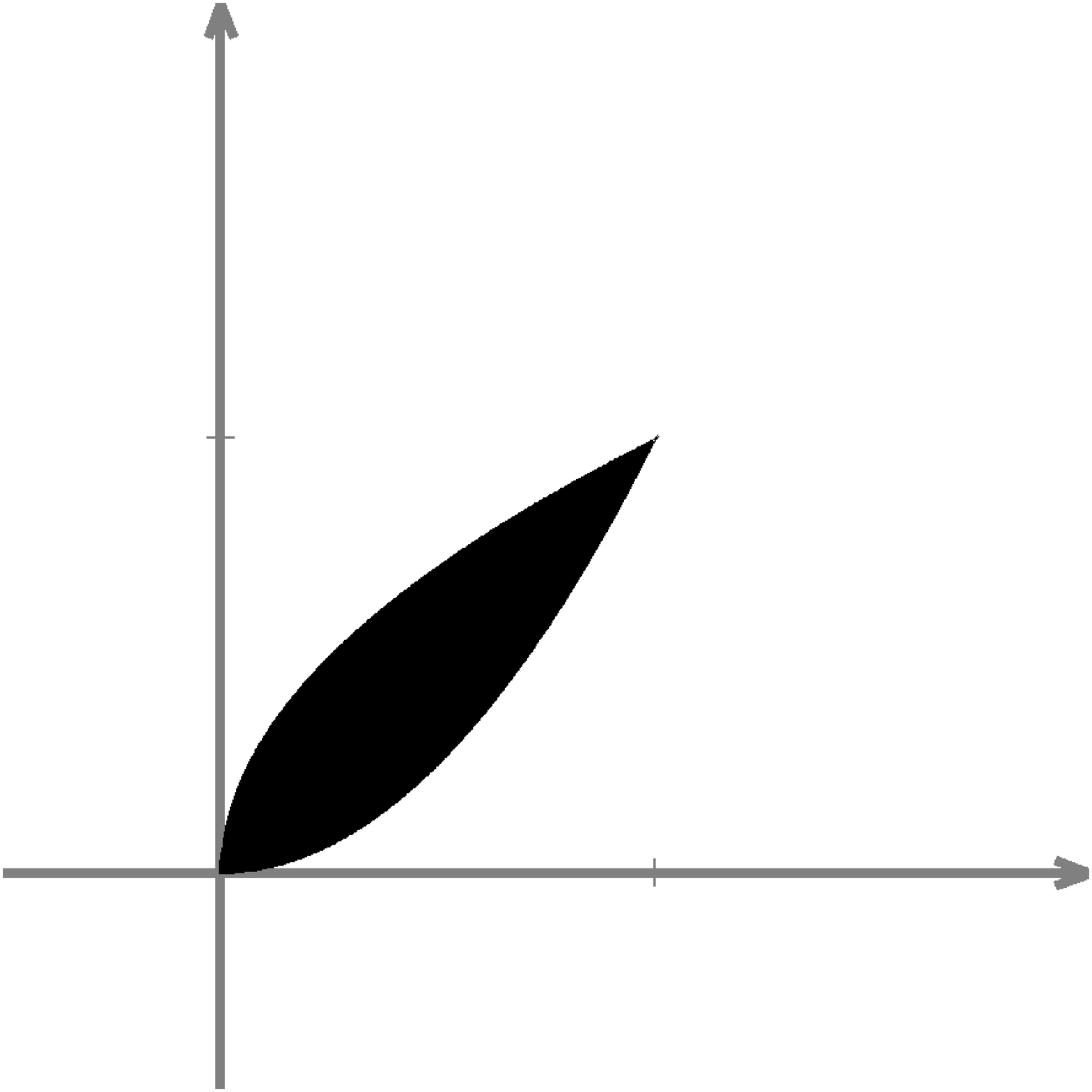}
\hspace{0.5cm}
\includegraphics[width=\figuresize]{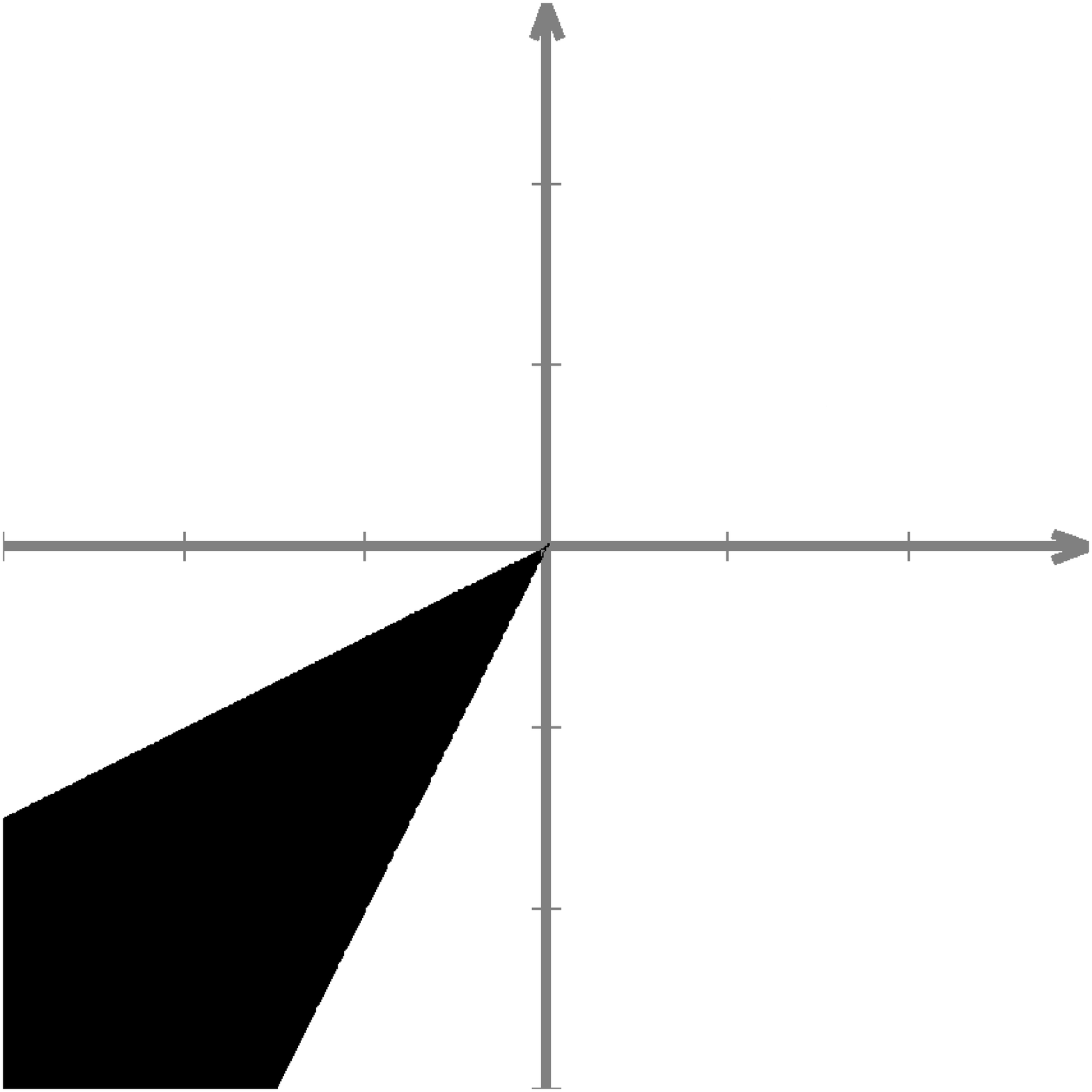}
\caption{$V = \{(x,y) \in {(\erre_{>0})}^2 \ |\ x^2 \leq y \leq \sqrt{x} \}$ (left picture), then $\ameba_0(V) = \{(x,y) \in \erre^2 \ |\ 2x \leq y \leq \frac{1}{2} x \}$ (right picture).} \label{fig:pow}
\end{center}
\end{figure}

Some examples of logarithmic limit sets are in figures \ref{fig:sin} and \ref{fig:pow}.

\begin{figure}[p]
\begin{center}
\includegraphics[height=\figuresize]{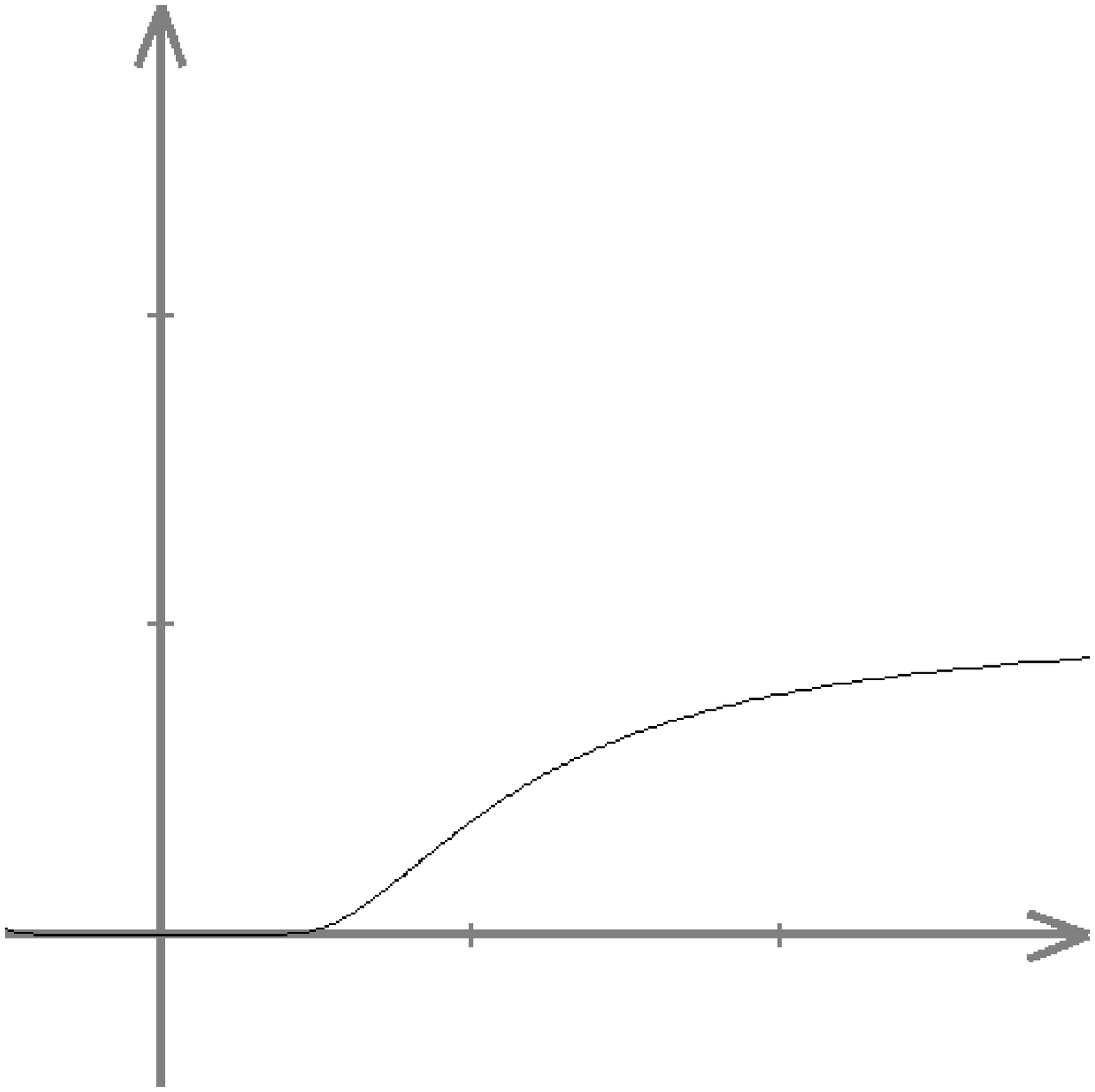}
\hspace{0.5cm}
\includegraphics[width=\figuresize]{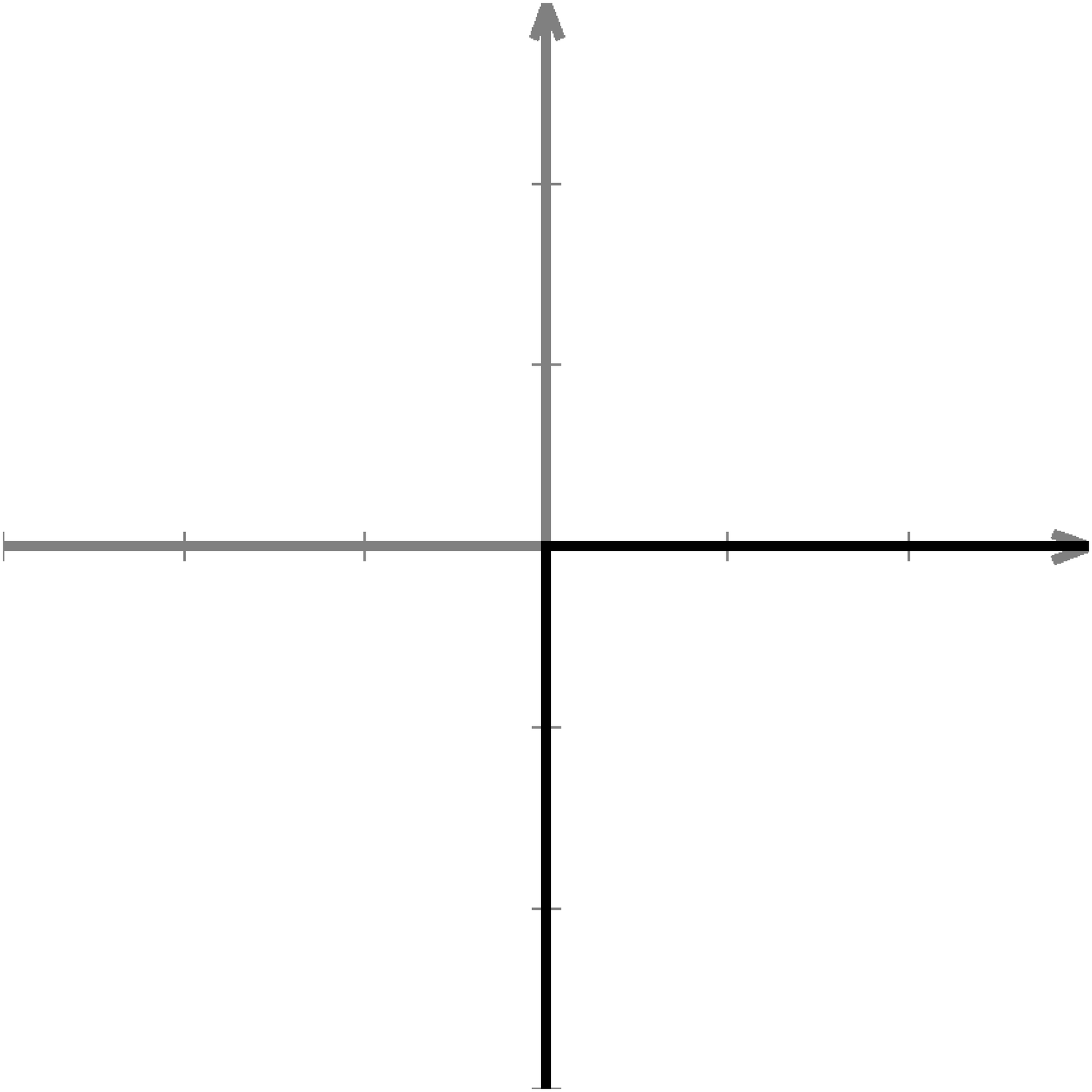}
\caption{$V = \{(x,y) \in {(\erre_{>0})}^2 \ |\ y = e^{-\frac{1}{x^2}} \}$ (left picture), then  $\ameba_0(V) = \{(x,y) \in \erre^2 \ |\ y=0, x \geq 0 \mbox{ or } x=0, y \leq 0 \}$ (right picture).}  \label{fig:exp}
\end{center}
\end{figure}

\begin{proposition} \label{prop:loglimsetproperties}
Given a set $V \subset {(\erre_{>0})}^n$ the following properties hold:
\begin{enumerate}
\item The logarithmic limit set $\ameba_0(V)$ is closed and $y \in \ameba_0(V)$ if and only if there exist a sequence $x(k)$ in $V$, and a sequence $t(k)$ in $(0,1)$ such that $t(k) \tende 0$ and $$\Log_{\left(\frac{1}{t(k)}\right)}(x(k))  \tende y$$
\item The logarithmic limit set $\ameba_0(V)$ is a cone in $\erre^n$.
\item We have that $0 \in \ameba_0(V)$ if and only if $V \neq \emptyset$. Moreover, $\ameba_0(V) =  \{0\}$ if and only if $V$ is compact and non-empty.
\item If $W \subset \erre^n$ we have $\ameba_0(V \cup W) = \ameba_0(V) \cup \ameba_0(W)$ and $\ameba_0(V \cap W) \subset \ameba_0(V) \cap \ameba_0(W)$. 
\end{enumerate}
\end{proposition}

\begin{proof} The first assertion is simply a restatement of proposition \ref{prop:generic properties}.

For the second one, we want to prove that if $\lambda > 0$ and $y \in \ameba_0(V)$, then $\lambda^{-1}y \in \ameba_0(V)$. There exists a sequence $x(k)$ in $V$ and a sequence $t(k)$ in $(0,1)$ such that $t(k) \tende 0$ and $\Log_{\left(\frac{1}{t(k)}\right)}(x(k)) \tende y$. Consider the sequence $x(k)$ and the sequence ${t(k)}^\lambda$. Now 
$$\Log_{\left(\frac{1}{{t(k)}^\lambda}\right)}(x(k)) = \dfrac{-1}{\log_e({t(k)}^\lambda)} \Log_e(x(k)) = \lambda^{-1}\dfrac{-1}{\log_e(t(k))} \Log_e(x(k))$$
and this sequence converges to $\lambda^{-1} y$.

The third and fourth assertions are trivial.
\end{proof}

Given a closed cone $C \subset \erre^n$, there is always a set $V \subset {(\erre_{>0})}^n$ such that $C = \ameba_0(V)$, simply take $V = \Log_e^{-1}(C)$. Then $\ameba_t(V) = C$ for all $t$.

Let $A = (a_{i j})\in GL_n(\erre)$. The matrix $A$ acts on $\erre^n$ in the natural way and, via conjugation with the map $\Log_e$, it acts on ${(\erre_{>0})}^n$. Explicitly, it induces the maps $A:\erre^n \freccia \erre^n$ and $\overline{A}:{(\erre_{>0})}^n \freccia {(\erre_{>0})}^n$:
$$A(x) = A(x_1, \dots, x_n) = (a_{11} x_1 + \dots + a_{1n} x_n,\ \dots \ , a_{n1} x_1 + \dots + a_{nn} x_n)$$
$$ \overline{A}(x) = \Exp_e \circ A \circ \Log_e(x) = 
(x_1^{a_{11}} x_2^{a_{12}} \cdots x_n^{a_{1n}}\ , \ \dots \ ,\ x_1^{a_{n1}} x_2^{a_{n2}} \cdots x_n^{a_{nn}} ) $$
If $V \subset {(\erre_{>0})}^n$ and $B \in GL_n(\erre)$, then $B(\ameba_0(V)) = \ameba_0(\overline{B}(V))$.

\begin{lemma}   \label{lemma:special point convergence}
$(0, \dots, 0,-1) \in \ameba_0(V)$ if and only if there exists a sequence $y(k)$ in $V$ such that $y_n(k) \freccia 0$ and

{$$\forall N \in \enne: \exists k_0 \in \enne: \forall k > k_0: \forall i \in \{1,\dots,n-1\}:$$ 
$$y_n(k) < {(y_i(k))}^N \mbox{ and }\  y_n(k) < {(y_i(k))}^{-N}$$}
\end{lemma}

\begin{proof} Suppose that $(0, \dots, 0,-1) \in \ameba_0(V)$, then by proposition \ref{prop:loglimsetproperties} there exists a sequence $y(k)$ in $V$ and a sequence $t(k)$ in  $(0,1)$ such that $t(k) \tende 0$ and $\Log_{\left(\frac{1}{t(k)}\right)}(y(k))  \tende (0, \dots, 0,-1)$. This means that 
$$\frac{-1}{\log_e(t(k))} \log_e(y_i(k)) \tende 
\left\{
\begin{array}{ll}
-1 & \mbox{ if } i = n\\
0  & \mbox{ if } i \in \{1,\dots,n-1\}
\end{array}
\right.$$
Now $t(k) \tende 0$ hence $\frac{-1}{\log_e(t(k))} \tende 0^+$, $\log_e(y_n(k)) \tende -\infty$, $y_n(k) \tende 0$ and
$$\forall i \in \{1,\dots,n-1\}: \frac{\log_e(y_i(k))}{\log_e(y_n(k))} \tende 0  $$ 
Hence $\forall \varepsilon > 0: \exists k_0: \forall k > k_0: \forall i < n:  {(y_n(k))}^\varepsilon < y_i(k) < {(y_n(k))}^{-\varepsilon} $.

We conclude by reversing the inequalities and choosing $\varepsilon = \frac{1}{N}$.

Conversely, if $y(k)$ has the stated property, then $|\Log_e(y(k))| \tende \infty$. It is possible to choose $t(k)$ such that $t(k) \tende 0$ and $\left|\Log_{\left(\frac{1}{t(k)}\right)}(y(k))\right| = 1$. Up to subsequences, the sequence $\Log_{\left(\frac{1}{t(k)}\right)}(y(k))$ converges to a point that, by reversing the calculations on first part of the proof, is $(0,\dots,0,-1)$. Hence $(0,\dots,0,-1) \in \ameba_0(V)$.
\end{proof}

\begin{corollary}
It follows that if there exists a sequence $x(k)$ in $V$ such that $x(k) \tende (a_1,\dots,a_{n-1},0)$, where $a_1, \dots, a_{n-1} > 0$, then $(0, \dots, 0,-1) \in \ameba_0(V)$. The converse is not true in general.
\end{corollary}

\begin{proof}
For the counterexample, see figure \ref{fig:exp}.
\end{proof}

We will see in theorem \ref{teo:loglimset} that if $V$ is definable in an o-minimal and polynomially bounded structure, the converse of the corollary becomes true.

A sequence $b(k)$ in ${(\erre_{>0})}^n$ is in \nuovo{standard position} in dimension $m$ if, denoted $g = n-m$, $b(k) \tende b = (b_1,\dots,b_g,0,\dots,0)$, with $b_1, \dots, b_g > 0$ and:
$$\forall N \in \enne: \exists k_0: \forall k > k_0: \forall i \in \{g+1, \dots, n-1\} : b_{i+1}(k) < {(b_i(k))}^N$$

\begin{lemma}   \label{lemma:coordinates}
Let $a(k)$ be a sequence in ${(\erre_{>0})}^n$ such that $a(k) \tende a = (a_1, \dots, a_h, 0,\dots,0)$, with $h < n$ and $a_1, \dots,a_h > 0$. There exists a subsequence (again denoted by $a(k)$) and a linear map $A:\erre^n \freccia \erre^n$ such that the sequence $b(k)=\left(\overline{A}(a(k))\right) \subset {(\erre_{>0})}^n$ is in standard position in dimension $m$, with $g = n-m \geq h$.
\end{lemma}

\begin{proof}
By induction on $n$. For $n = 1$ the statement is trivial. Suppose that the statement holds for $n-1$. Consider the logarithmic image of the sequence: $\Log(a(k))$. Up to extracting a subsequence, the sequence $\left(\frac{\Log(a(k))}{|\Log(a(k))|}\right)$ converges to a unit vector $v = (0,\dots,0,v_{h+1},\dots,v_n)$. There exists a linear map $B$, acting only on the last $n-h$ coordinates, sending $v$ to $(0,\dots,0,-1)$. By lemma \ref{lemma:special point convergence}, the map $\overline{B}$ sends $a(k)$ to a sequence $b(k)$ such that $b_n(k) \tende 0$ and $\forall N \in \enne: \exists k_0: \forall k > k_0: \forall i \in \{1,\dots, n-1\}: $
$$b_n(k) < (b_i(k))^N \mbox{ and } b_n(k) < (b_i(k))^{-N} $$
As $B$ only acted on the last $n-h$ coordinates, for $i \in \{1,\dots,h\}$, $b_i(k) \tende a_i \neq 0$. Up to subsequences we can suppose that for every $i \in \{h+1,\dots,n-1\}$ one of the three possibilities occur: $b_i(k) \tende 0$, $b_i(k) \tende b_i \neq 0$, $b_i(k) \tende +\infty$. Up to a change of coordinates with maps of the form 
$$B_i(x_1,\dots,x_i, \dots, x_n) = (x_1,\dots, -x_i, \dots, x_n) $$
$$\overline{B_i}(x_1,\dots,x_i, \dots, x_n) = (x_1,\dots,x_i^{-1}, \dots, x_n) $$
we can suppose that for every $i \in \{1,\dots,n-1\}$ either $b_i(k) \tende 0$ or $b_i(k) \tende b_i \neq 0$. Up to reordering the coordinates, we can suppose that exists $g \geq h$ such that for $i \in \{1,\dots,g\}$: $b_i(k) \tende b_i \neq 0$ and for $i > g$, $b_i(k) \tende 0$.
Now consider the projection on the first $n-1$ coordinates: $\pi:\erre^n \freccia \erre^{n-1}$. By inductive hypothesis there exists a linear map $C: \erre^{n-1} \freccia \erre^{n-1}$ sending the sequence $\pi(b(k))$ in a sequence $c(k)$ satisfying:
$$\forall N \in \enne: \exists k_0: \forall k > k_0: \forall i \in \{g+1, \dots, n-2\} : b_{i+1}(k) < {(b_i(k))}^N$$

The composition of $B$ and a map that preserves the last coordinate and acts as $C$ on the first ones is the searched map.
\end{proof}

The \nuovo{basic cone} defined by the vector $N = (N_1, \dots, N_{n-1}) \in \enne^{n-1}$ is:
$$B_N = \{ x \in \erre^n \ |\ \forall i: x_i \leq 0 \mbox{ and } \forall i < n : x_{i+1} \leq N_i x_i \} $$

Note that if $N' = (N_1',\dots, N_{n-1}')$, with $\forall i: N_i' \geq N_i$, then $B_N' \subset B_N$.

The \nuovo{exponential basic cone} in ${(\erre_{>0})}^n$ defined by the vector $N = (N_1, \dots, N_{n-1}) \in \enne^{n-1}$ and the scalar $h>0$ is the set:
$$E_{N,h} = \{ x \in \erre^n \ |\ \forall i: 0 < x_i \leq h \mbox{ and } \forall i < n : x_{i+1} \leq {x_i}^{N_i} \}$$

\begin{lemma} The following easy facts about basic cones holds:
\begin{enumerate}
\item The logarithmic limit set of an exponential basic cone is a basic cone:
$$\ameba_0( E_{N,h}) = B_N $$
\item If $b(k) \subset {(\erre_{>0})}^n$ is in standard position in dimension $n$, and $E_{N,h}$ is an exponential basic cone, then for large enough $k$, $b(k) \in E_{N,h}$.
\end{enumerate}
\end{lemma}

\subsection{Definable sets in o-minimal structures}  \label{subsez:o-minimal structures}

In this subsection we report some notations ans some definitions of model theory and o-minimal geometry we will use later, see \cite{EFT84} and \cite{Dr} for details. Given a \nuovo{set of symbols} $S$ (see \cite[chap. II, def. 2.1]{EFT84}), we denote by $L_S$ the corresponding \nuovo{first order language} (see \cite[chap. II, def. 3.2]{EFT84}). If $S'$ is an \nuovo{expansion} of a set of symbols $S$ we will write $S \subset S'$. The theory of real closed fields uses the set of symbols of ordered semirings: $\mathcal{OS} = (\{\leq\},\{+,\cdot\}, \emptyset)$ or, equivalently, the set of symbols of ordered rings $\mathcal{OR} = (\{\leq\},\{+,-,\cdot\},\{0,1\})$, an expansion of $\mathcal{OS}$. In the following we will use these sets of symbols and some of their expansions.

We usually will denote an $S$-structure by $\overline{M} = (M,a)$, where $M$ is a set, and $a$ is the  \nuovo{interpretation}, (see \cite[cap. III, def. 1.1]{EFT84}). Given an $S$-structure $\overline{M} = (M,a)$, and an $L_S$-formula $\phi$ without free variables, we will write $\overline{M} \vDash \phi$ if $\overline{M}$ satisfies $\phi$ (see \cite[chap. III, def. 3.1]{EFT84}).

A real closed field can be defined as an $\mathcal{OS}$- or an $\mathcal{OR}$-structure satisfying a suitable infinite set of first order axioms. The natural $\mathcal{OS}$-structure on $\erre$ will be denoted by $\overline{\erre}$.

If $\overline{M} = (M,a)$ is an $S$-structure, and $S'$ is an expansion of $S$, an $S'$-structure $(M,a')$ is an \nuovo{expansion} of the $S$-structure $(M,a)$ if $a'$ restricted to the symbols of $S$ is equal to $a$. If $M \subset N$, an $S$-structure $\overline{N} = (N,b)$ is an \nuovo{extension} of an $S$-structure $\overline{M} = (M,a)$ if for all $ s \in S$, $b(s)_{|M} = a(s)$.

A \nuovo{definable subset} of $M^n$ is a set that is defined by an $(L_S)$-formula $\phi(x_1,\dots,x_n,y_1,\dots,y_m)$ and by parameters $a_1, \dots a_m \in M$, and a \nuovo{definable map} is a map whose graph is definable. For example if $M$ is an $\mathcal{OS}$-structure satisfying the axioms of real closed fields, the definable sets are the semi-algebraic sets, and the definable maps are the semi-algebraic maps.

Let $S$ be an expansion of $\mathcal{OS}$, and let $\mathfrak{R} = (\erre,a)$ be an $S$-structure that is an \nuovo{o-minimal} and \nuovo{polynomially bounded} expansion of $\overline{\erre}$. In \cite{Mi94} it is shown that if $f:\erre \freccia \erre$ is definable and not ultimately $0$, there exist $r,c \in \erre$, $c \neq 0$, such that
$$\displaystyle \lim_{x \tende +\infty} \dfrac{f(x)}{x^r} = c $$
The set of all such $r$ is a subfield of $\erre$, called the \nuovo{field of exponents} of $\mathfrak{R}$. For example the $\mathcal{OR}$-structure $\overline{\erre}$ is polynomially bounded with field of exponents $\qu$.

If $\Lambda \subset \erre$ is a subfield, we can construct an expansion of $S$ and $\mathfrak{R}$ by adding the power functions with exponents in $\Lambda$. We expand $S$ to $S^{\Lambda}$ by adding a function symbol $f_\lambda$ for every $\lambda \in \Lambda$, and we expand $\mathfrak{R}$ to an $S^{\Lambda}$-structure $\mathfrak{R}^{\Lambda}$ interpreting the function symbol $f_\lambda$ by the function that is $x \freccia x^{\lambda}$ for positive numbers and $x \freccia 0$ on negative ones. The structure $\mathfrak{R}^{\Lambda}$ is again o-minimal, as its definable sets are definable in the structure $(\mathfrak{R},e^x)$, that is o-minimal by \cite{Spe99}.

Suppose that the expansion of $\mathfrak{R}$ constructed by adding the family of functions ${\{x^r_{|[1,2]}\}}_{r \in \Lambda}$ is polynomially bounded, then $\mathfrak{R}^{\Lambda}$ is too (see \cite{Mi03}). For example if $(\mathfrak{R},e^x_{|[0,1]})$ is polynomially bounded, then $\mathfrak{R}^{\Lambda}$ is too.

In the following we will work with o-minimal, polynomially bounded structures $\mathfrak{R}$ expanding $\overline{\erre}$, with the property that $\mathfrak{R}^\erre$ is polynomially bounded. We will call such structures \nuovo{regular polynomially bounded structures}.

One example of regular polynomially bounded structure is $\erre_{\mathrm an}$, the real numbers with restricted analytic functions, see \cite{DMM94} for details. This structure has field of exponents $\qu$, while $\erre_{\mathrm an}^{\erre}$ has field of exponents $\erre$.

As $\erre_{\mathrm an}$ is an expansion of $\overline{\erre}$, also $\overline{\erre}^\erre$ is polynomially bounded, hence $\overline{\erre}$ is a regular polynomially bounded structure.

Other regular polynomially bounded structures we will not use here are the structure $\erre_{\mathrm an^*}$ of the real field with convergent generalized power series, (see \cite{DS98}), the field of real numbers with multisummable series (see \cite{DS00}) and the structures defined by a quasianalytic Denjoy-Carlemann class (see \cite{RSW03}).

\section{Logarithmic limit sets of definable sets}           \label{sez:poly}

\subsection{Some properties of definable sets}

Let $\mathfrak{R}$ be an o-minimal and polynomially bounded expansion of $\overline{\erre}$.

\begin{lemma}
For every definable function $f:{(\erre_{>0})}^n \freccia \erre_{>0}$, there is a basic exponential cone $C$ and $N \in \enne$ such that $f_{|C}(x_1, \dots, x_n) \geq {(x_n)}^N $.
\end{lemma}

\begin{proof}
Fix a basic exponential cone $C \subset {(\erre_{>0})}^n$. By the \L{}ojasiewicz inequality (see \cite[4.14]{DM96}) there exist $N \in \enne$ and $Q > 0$ such that $Q f_{|C}(x_1, \dots, x_n) \geq {(x_n)}^N $. The thesis follows by choosing an exponent bigger than $N$ and a suitable basic exponential cone smaller than $C$.
\end{proof}

\begin{lemma} \label{lemma:celldecomposition}
Every cell decomposition of ${(\erre_{>0})}^n$ has a cell containing a basic exponential cone. 
\end{lemma}

\begin{proof} This proof is based on the cell decomposition theorem, see \cite[chap. 3]{Dr} for details. By induction on $n$. For $n=1$, the statement is trivial.

Suppose the lemma true for $n$. If $\{C_i\}$ is a cell decomposition of ${(\erre_{>0})}^{n+1}$, and if $\pi:{(\erre_{>0})}^{n+1} \freccia {(\erre_{>0})}^n$ is the projection on the first $n$ coordinates, then $\{\pi(C_i)\}$ is a cell decomposition of ${(\erre_{>0})}^n$, hence, by induction, it contains a basic exponential cone $D$ of ${(\erre_{>0})}^n$. Then $\pi^{-1}(D) \times (0,1]$ contains a cell of the form
$$E = \{(\bar{x},x_{n+1}) \ |\ \bar{x} \in D, 0 < x_{n+1} < f(\bar{x}) \} $$
where $\bar{x}=(x_1, \dots, x_n)$ and $f:D \freccia (0,1]$ is definable. By previous lemma, there is a basic exponential cone $D' \subset D$ and $N \in \enne$ such that $f_{|D'}(\bar{x}) \geq {(x_n)}^N$. Hence $E$ contains the basic exponential cone
$$\{(\bar{x},x_{n+1}) \ |\ \bar{x} \in D', 0 < x_{n+1} \leq {(x_n)}^N \}$$
\end{proof}

\begin{corollary}  \label{corol:almostnotzero}
Let $V \subset {(\erre_{>0})}^n$ be definable in $\mathfrak{R}$, and suppose that $V$ contains a sequence $x(k)$ in standard position in dimension $n$. Then $V$ contains an exponential cone.
\end{corollary}

\begin{proof}
Let $\{C_i\}$ be a cell decomposition of ${(\erre_{>0})}^n$ adapted to $V$. By previous lemma, one of the cells contains an exponential cone $D$. By  hypothesis, if $k$ is sufficiently large, $x(k) \in D$, hence $D \subset V$.
\end{proof}

\begin{corollary}  \label{corol:almostnotzeroplane}
Let $V \subset {(\erre_{>0})}^2$ be definable in $\mathfrak{R}$, and suppose that exists a sequence $x(k)$ in $V$ such that $x(k) \tende 0$ and 
$$\forall N \in \enne: \exists k_0: \forall k > k_0: x_2(k) < {(x_1(k))}^N$$ 
Then there exist $h_0 > 0$ and $M \in \enne$ such that 
$$\{ x \in \erre^2 \ |\ 0 < x_1 < h_0 \mbox{ and } 0 < x_2 < {(x_1)}^M \} \subset V$$ 
\end{corollary}

\begin{proof}
This is precisely the previous corollary with $n = 2$.
\end{proof}

\begin{lemma}   \label{lemma:almostnotzero}
Let $V \subset {(\erre_{>0})}^n$ be definable in $\mathfrak{R}$, and suppose that there exists a sequence $x(k)$ in $V$, an integer $m \in \{1,\dots, n\}$ and, denoted $g=n-m$, positive numbers $a_1, \dots, a_g > 0$ such that $x(k) \tende (a_1,\dots, a_g, 0, \dots, 0)$, and such that:
$$\forall N \in \enne: \exists k_0: \forall k > k_0: \forall i \in \{ g+1, \dots, n-1\} : x_n(k) < {(x_i(k))}^N$$ 
Then for every $\varepsilon > 0$ there exist a sequence $y(k)$ in $V$ and positive real numbers $b_1, \dots b_{n-1} > 0$ such that $y(k) \tende (b_1, \dots b_{n-1},0)$ and for all $i \in \{1, \dots g\}$ we have $|b_i-a_i| < \varepsilon$.
\end{lemma}

\begin{proof}
If $n = 2$ the statement follows by corollary \ref{corol:almostnotzeroplane}. By induction on $n$ we suppose the statement true for definable sets in $\erre^{n'}$ with $n' < n$. We split the proof in two cases, when $m < n$ and when $m=n$.

If $m < n$, fix an $\varepsilon > 0$, smaller than every $a_i$, and consider the parallelepiped
$$c_\varepsilon = \{ (z_1,\dots,z_n) \in \erre^n \ |\  |z_1-a_1| < {\textstyle\frac{1}{2}}\varepsilon, \dots, |z_g-a_g| < {\textstyle\frac{1}{2}}\varepsilon \} $$
Let $\pi:\erre^n \freccia \erre^{n-m}$ be the projection on the last $n-m$ coordinates. The set $\pi(V \cap c_\varepsilon)$ is definable in $\erre^m$, the sequence $\pi(x(k))$ satisfies the hypotheses of the lemma, hence, by induction, there exists a sequence $z(k) \in \pi(V \cap c_\varepsilon)$ converging to the point $(b_{g+1},\dots,b_{n-1},0)$. Let $y(k)$ be a sequence such that $y(k) \in \pi^{-1}(z(k))$. We can extract a subsequence (called again $y(k)$) such that $y(k) \tende (0,b_2, \dots b_n)$ where for all $i \in \{1, \dots g\}$ we have $|b_i-a_i| \leq \frac{1}{2}\varepsilon$.

If $m=n$, then $x(k) \tende 0$. The sequence $\left(\frac{(x_1(k),\dots,x_{n-1}(k))}{|(x_1(k),\dots,x_{n-1}(k))|}\right)$ is contained in the unit sphere $S^{n-2}$, and, up to subsequences, we can suppose that it converges to a unit vector $v=(v_1,\dots,v_{n-1}) \in {(\erre_{\geq 0})}^{n-1}$. Up to reordering, $v = (v_1, \dots, v_h, 0,\dots, 0)$, with $v_1, \dots, v_h > 0$. Fix an $\alpha > 0$ and consider the cone 
$$C_v(\alpha) = \{ y \in \erre^n \ |\ \frac{\left<(y_1, \dots, y_{n-1}),v\right>}{|(y_1, \dots, y_{n-1})||v|} > \cos \alpha \} $$
Let $\pi:\erre^n \freccia \erre^{n-h}$ the projection on the last $n-h$ coordinates. The set $\pi(V \cap C_v(\alpha))$ is definable in $\erre^{n-h}$, the sequence $\pi(x(k))$ satisfies the hypotheses of the lemma, hence, by induction, there exists a sequence $z(k) \in \pi(V \cap C_v(\alpha))$ converging to the point $(b_{h+1},\dots,b_{n-1},0)$. Let $y(k)$ be a sequence such that $y(k) \in \pi^{-1}(z(k))$. Up to subsequences, $y(k) \tende (b_1, \dots b_{n-1},0)$. As $y(k) \in C_v(\alpha)$, for all $i > h$, if $y_i(k) \tende 0$, then $y(k) \tende 0$. As $b_{h+1}, \dots, b_{n-1} > 0$, then also $b_1, \dots, b_h > 0$.
\end{proof}

\subsection{Polyhedral structure}      \label{subsez:polyh struc}

Let $V \subset {(\erre_{>0})}^n$ be a definable set. Our main object of study is $\ameba_0(V)$, the logarithmic limit set of $V$. Suppose that $\mathfrak{R}$ has field of exponents $\Lambda \subset \erre$. Given a matrix , the set $B(\ameba_0(V))$ is the logarithmic limit set of $\overline{B}(V)$. The components of $\overline{B^{-1}}$ are all definable in $\mathfrak{R}$ because their exponents are in $\Lambda$, hence the set $\overline{B}(V)$ is again definable.

\begin{theorem}  \label{teo:loglimset}
Let $V \subset ({\erre_{>0}})^n$ be a set definable in an o-minimal and polynomially bounded structure. The point $(0, \dots, 0,-1)$ is in $\ameba_0(V)$ if and only if there exists a sequence $x(k)$ in $V$ such that $x(k) \tende (a_1,\dots,a_{n-1},0)$, where $a_1, \dots, a_{n-1} > 0$.
\end{theorem}

\begin{proof}
If there exists such an $x(k)$, then it is obvious that $(0, \dots, 0,-1) \in \ameba_0(V)$. Vice versa, if $(-1,0, \dots, 0) \in \ameba_0(V)$, then by lemma \ref{lemma:special point convergence} there exists a sequence $y(k)$ in $V$ such that $y_n(k) \freccia 0$ and

{$$\forall N \in \enne: \exists k_0 \in \enne: \forall k > k_0: \forall i \in \{1,\dots,n-1\}:$$ 
$$y_n(k) < {(y_i(k))}^N \mbox{ and } y_n(k) < {(y_i(k))}^{-N}$$}
Up to subsequences we can suppose that for all $i \in \{1,\dots,{n-1}\}$ one of the three possibilities occur: $y_i(k) \tende 0$, $y_i(k) \tende a_i \neq 0$, $y_i(k) \tende +\infty$. Up to a change of coordinates with maps of the form 
$$B_i(x_1,\dots,x_i, \dots, x_n) = (x_1,\dots, -x_i, \dots, x_n) $$
we can suppose that for all $i \in \{1,\dots,n-1\}$ either $y_i(k) \tende 0$ or $y_i(k) \tende a_i \neq 0$. Then we can apply lemma \ref{lemma:almostnotzero}, and we are done.
\end{proof}

Now we suppose that $\mathfrak{R}$ is a regular polynomially bounded structure, or, equivalently, that $\mathfrak{R}$ has field of exponents $\erre$. Let $x \in \ameba_0(V)$. We want to describe a neighborhood of $x$ in $\ameba_0(V)$. To do this, we choose a map $B \in GL_n(\erre)$ such that $B(x) = (0, \dots, 0,-1)$. Now we only need to describe a neighborhood of $(0, \dots, 0,-1)$ in $\ameba_0(\overline{B}(V))$. As logarithmic limit sets are cones, we only need to describe a neighborhood of $0$ in 
$$ H = \{(x_1, \dots, x_{n-1}) \in \erre^{n-1} \ |\ (x_1, \dots, x_{n-1},-1) \in \ameba_0(\overline{B}(V)) \}$$
We define a one-parameter family and its limit:
$$W_t = \{ (x_1, \dots, x_{n-1}) \in {(\erre_{>0})}^{n-1} \ |\ (x_1, \dots, x_{n-1},t) \in \overline{B}(V) \} $$
$$W = \left( \lim_{t \tende 0} W_t \right) \cap {(\erre_{>0})}^{n-1}$$
The set $W$ is a definable subset of ${(\erre_{>0})}^{n-1}$. Its logarithmic limit set is denoted, as usual, by $\ameba_0(W)\subset \erre^{n-1}$. By previous theorem, as $(0, \dots, 0,-1)\in\ameba_0(\overline{B}(V))$, $W$ is not empty, hence $0\in \ameba_0(W)$. We want to prove that there exists a neighborhood $U$ of $0$ in $\erre^{n-1}$ such that $\ameba_0(W) \cap U = H \cap U$ or, in other words, that $\ameba_0(W) \cap H$ is a neighborhood of $0$ both in $\ameba_0(W)$ and $H$.

A \nuovo{flag} in $\erre^n$ is a sequence $(V_0, V_1, \dots, V_h)$, $h \leq n$, of subspaces of $\erre^n$ such that $V_0 \subset V_1 \subset \dots \subset V_h \subset \erre^n$ and $\dim V_i = i$. We say that a sequence $x(k)$ in $\erre^n$ converges to the point $y$ \nuovo{along} the flag $(V_1, V_2, \dots, V_h)$ if $x(k) \tende y$, $\forall k: x(k) - y \in V_h \setminus V_{h-1}$ and $\forall i \in \{0, \dots, h-2\}$, the sequence $\pi_i(x(k))$ converges to the point $\pi_i(V_{i+1})$, where $\pi_i:V_h\setminus V_i \freccia \pro(V_h / V_i)$ is the canonical projection.

\begin{lemma} \label{lemma:flag convergence}
For all sequences $x(k)$ in $\erre^n$ converging to a point $y$, there exists a flag $(V_0, \dots, V_h)$ and a subsequence of $x(k)$ converging to $y$ along $(V_0, \dots, V_h)$.
\end{lemma}

\begin{proof}
It follows from the compactness of $\pro(V_h / V_i)$.
\end{proof}

\begin{lemma}
Let $x(k) \subset H$ be a sequence converging to $0$. Then at least one of its points is in $\ameba_0(W) \cap H$.
\end{lemma}

\begin{proof}
By lemma \ref{lemma:flag convergence}, we can extract a subsequence, again denoted by $x(k)$, converging to zero along a flag $(V_0,V_1, \dots, V_{h})$ in $\erre^{n-1}$. Up to a linear change of coordinates, we can suppose that this flag is given by $(\{0\}, \Span(e_{n-1}), \Span(e_{n-2},e_{n-1}), \dots, \Span(e_{n-h}, \dots, e_{n-1}))$. Hence for $i \in \{1,\dots,n-h-1\}$ we have $x_i(k) = 0$. Again by extracting a subsequence and by a change of coordinates with maps of the form 
$$B_i(x_1,\dots,x_i, \dots, x_n) = (x_1,\dots, -x_i, \dots, x_n) $$ 
with $i \in \{n-h,\dots,n-1\}$, we can suppose that for all such $i$, $x_i(k) < 0$.

By proposition \ref{prop:loglimsetproperties}, as $H \subset \ameba_0(B(V))$, for every point $x(k)$ there exists a sequence $y(k,l)$ in $B(V)$ and a sequence $t(k,l)$ in $(0,1)$ such that $t(k,l) \tende 0$ and $\Log_{\left(\frac{1}{t(k,l)}\right)}(y(k,l))  \tende x(k)$. By theorem \ref{teo:loglimset} we can choose $y(k,l)$ such that $y(k,l) \tende a(k)$, with $a_i(k) > 0$ for $i \in \{1,\dots,n-h-1\}$, and $a_i(k) = 0$ for $i \in \{n-h,\dots,n\}$. Up to a change of coordinates with maps of the form 
$$B_i(x_1,\dots,x_i, \dots, x_n) = (x_1,\dots, -x_i, \dots, x_n) $$
with $i \in \{1,\dots,n-h-1\}$, we can suppose that the sequence $a(k)$ is bounded and that, up to subsequences, it converges to a point $a$, with $a_i = 0$ for $i \in \{n-h,\dots,n\}$.

Let $\pi:\erre^{n} \freccia \erre^{n-h-1}$ be the projection on the first $n-h-1$ coordinates. Then $\pi(a(k)) \subset {(\erre_{>0})}^{n-h-1}$. By lemma \ref{lemma:coordinates} we can suppose that $\pi(a(k))$ is in standard position, i.e. $a_1, \dots, a_g > 0, a_{g+1} = \dots = a_{n} = 0$ and:
$$\forall N \in \enne: \exists k_0: \forall k > k_0: \forall i \in \{g+1, \dots, n-h-2\} : a_{i+1}(k) < {(a_i(k))}^N$$

From the sequences $y(k,l)$, we extract a diagonal subsequence $z(k)$ in the following way. For every $k$, the sequence $y(k,l)$ converges to $a(k) = (a_1(k), \dots, a_{n-h-1}(k), 0,\dots,0)$. As $\Log_{\frac{1}{t(k,l)}}(y(k,l))  \tende x(k) = (x_1(k),\dots,x_{n-h-1}(k),0,\dots,0,-1)$, for all $i \in \{n-h,\dots,n\}$ we have
$$\frac{\log(y_1(k,l))}{\log(y_i(k,l))} \freccia \frac{-1}{x_i(k)}$$
We can choose an $l_0$ such that:
\begin{enumerate}
\item $\forall i \in \{1,\dots,h\}: \left|\frac{\log(y_1(k,l_0))}{\log(y_i(k,l_0))} - \frac{-1}{x_i(k)}\right| < \frac{1}{k}$
\item $|y(k,l_0) - a(k)| < \frac{1}{k} $
\end{enumerate}
We define $z(k) = y(k,l_0)$. Now $z(k) \tende a=(a_1, \dots, a_g, 0,\dots,0)$ and, as $x(k) \tende (0,\dots,0,-1)$ along the flag $(V_0, \dots, V_h)$, we have:
$$\forall N \in \enne: \exists k_0: \forall k > k_0: \forall i \in \{g+1, \dots, n-1\} : x_{i+1}(k) < {(x_i(k))}^N$$

Let $r$ be smaller than every $a_1,\dots,a_g$. Consider the parallelepiped
$$c_r = \{ (z_1,\dots,z_n) \in \erre^n \ |\  |z_1-a_1| < {\textstyle\frac{1}{2}}r, \dots, |z_g-a_g| < {\textstyle\frac{1}{2}}r \}  $$
Let $\pi:\erre^n \freccia \erre^{n-g}$ be the projection on the last $n-g$ coordinates. The set $\pi(B(V) \cap c_r)$ is definable in $\erre^{n-g}$, and the sequence $\pi(z(k))$ satisfies the hypotheses of corollary \ref{corol:almostnotzero}, hence $\pi(B(V) \cap c_r)$ contains a basic exponential cone, hence $\pi(W \cap c_r)$ also contains one. This means that $\ameba_0( \pi(W) \cap c_r )$ contains a basic cone. Hence also $\ameba_0(W \cap c_r)$ contains this cone, and also $\ameba_0(W)$. At least one of the points $x(k)$ is in this cone.
\end{proof}

\begin{lemma}  \label{lemma:ray of coincidence}
Let $x \in \ameba_0(W)$. Then the number 
$$r(x) = \sup\{ r \ |\ \forall 0 \leq \lambda \leq r: \lambda x \in H \} $$
is strictly positive.
\end{lemma}

\begin{proof}
Let $x \in \ameba_0(W)$. By a linear change of coordinates, we can suppose that $x=(0,\dots,0,-1,-1)$. By theorem \ref{teo:loglimset} there is a sequence $x(k)$ in $W$ converging to the point $(a_1,\dots,a_{n-2},0)$, with $a_1,\dots,a_{n-2} > 0$. As $W$ is the limit of the family $W_t$, for every $k$ there is a sequence $y(k,l)$ in $B(V)$ converging to $(x(k),0)$. We can construct a diagonal sequence $z(k)$ in the following way: for every $k$ we can choose an $l_0$ such that
$$|y(k,l_0) - (x(k),0)| < {(x_{n-1}(k))}^k $$
The sequence $z(k)$ converges to $(a_1,\dots,a_{n-2},0,0)$. Let $r$ be smaller that any of the $a_1,\dots,a_{n-2}$. Consider the parallelepiped
$$c_r = \{ (z_1,\dots,z_n) \in \erre^n \ |\  |z_1-a_1| < {\textstyle\frac{1}{2}}r, \dots, |z_{n-2}-a_{n-2}| < {\textstyle\frac{1}{2}}r \}  $$
Let $\pi:\erre^n \freccia \erre^2$ be the projection on the last $2$ coordinates. The set $\pi(B(V) \cap c_r)$ is definable in $\erre^2$, and the sequence $\pi(z(k))$ satisfies the hypotheses of corollary \ref{corol:almostnotzeroplane}, hence $\pi(B(V) \cap c_r)$ contains a basic exponential cone. This means that there exists a number $r'>0$ such that
$$\{(0,\dots,0,z,-1) \ |\ -r' \leq z \leq 0 \} \subset H$$
\end{proof}

\begin{theorem}   \label{teo:neighborhood}
Let $V \subset ({\erre_{>0}})^n$ be a set definable in a regular polynomially bounded structure. Let $x \in \ameba_0(V)$ and choose a map $B \in GL_n(\erre)$ such that $B(x) = (0, \dots, 0,-1)$. We recall that
$$ H = \{(x_1, \dots, x_{n-1}) \in \erre^{n-1} \ |\ (x_1, \dots, x_{n-1},-1) \in \ameba_0(\overline{B}(V)) \}$$
$$W_t = \{ (x_1, \dots, x_{n-1}) \in {(\erre_{>0})}^{n-1} \ |\ (x_1, \dots, x_{n-1},t) \in \overline{B}(V) \} $$ 
$$W = \left( \lim_{t \tende 0} W_t \right) \cap {(\erre_{>0})}^{n-1}$$ 
Then there exists a neighborhood $U$ of $0$ in $\erre^{n-1}$ such that $\ameba_0(W) \cap U = H \cap U$.
\end{theorem}

\begin{proof}
We will prove that $\ameba_0(W) \cap H$ is a neighborhood of $0$ both in $\ameba_0(W)$ and in $H$. Previous lemma implies that if $x(k)$ is a sequence in $H$ converging to $0$, then at least one of its points is in $\ameba_0(W)$, hence $\ameba_0(W) \cap H$ is a neighborhood of $0$ in $H$.

To prove that $\ameba_0(W) \cap H$ is also a neighborhood of $0$ in $\ameba_0(W)$, we only need to prove that if $r$ is the function defined in lemma \ref{lemma:ray of coincidence}, there exists an $\varepsilon > 0$ such that 
$$\forall x \in \ameba_0(W) \cap S^{n-2}: r(x) > \varepsilon $$
But this is true, because we already know that $\ameba_0(W) \cap H$ is a neighborhood of $0$ in $H$.
\end{proof}

\begin{theorem}   \label{teo:polyh complex}
Let $V \subset ({\erre_{>0}})^n$ be a set definable in a regular polynomially bounded structure. The logarithmic limit set $\ameba_0(V)$ is a polyhedral complex. Moreover, if $\dim V = m$, then $\dim \ameba_0(V) \leq m$.
\end{theorem}

\begin{proof}
By induction on $n$. For $n = 1$ the statement is trivial, as a cone in $\erre$ is a polyhedral set, and every zero dimensional definable set is compact, hence its logarithmic limit set is a point. Suppose the statement true for $n-1$. For every $x \in \ameba_0(V)$ there is a linear map $B$ sending $x$ to $(0,\dots,0,-1)$. The statement in \cite[4.7]{DM96} implies that the definable set $W \subset ({\erre_{>0}})^n$ has dimension less than or equal to $m-1$, hence $\ameba_0(W)$ is a polyhedral set of dimension less than or equal to $m-1$ (by inductive hypothesis). By previous theorem a neighborhood of the ray $\{\lambda x \ |\ \lambda \geq 0 \}$ in $\ameba_0(V)$ is the cone over a neighborhood of $0$ in $\ameba_0(W)$, hence it is a polyhedral complex of dimension less than or equal to $m$. By compactness of the sphere $S^{n-1}$, $\ameba_0(V)$ can be covered by a finite number of such neighborhoods, hence it is a polyhedral complex of dimension less than or equal to $m$.
\end{proof}

\begin{figure}[p]
\begin{center}
\includegraphics[height=\figuresize]{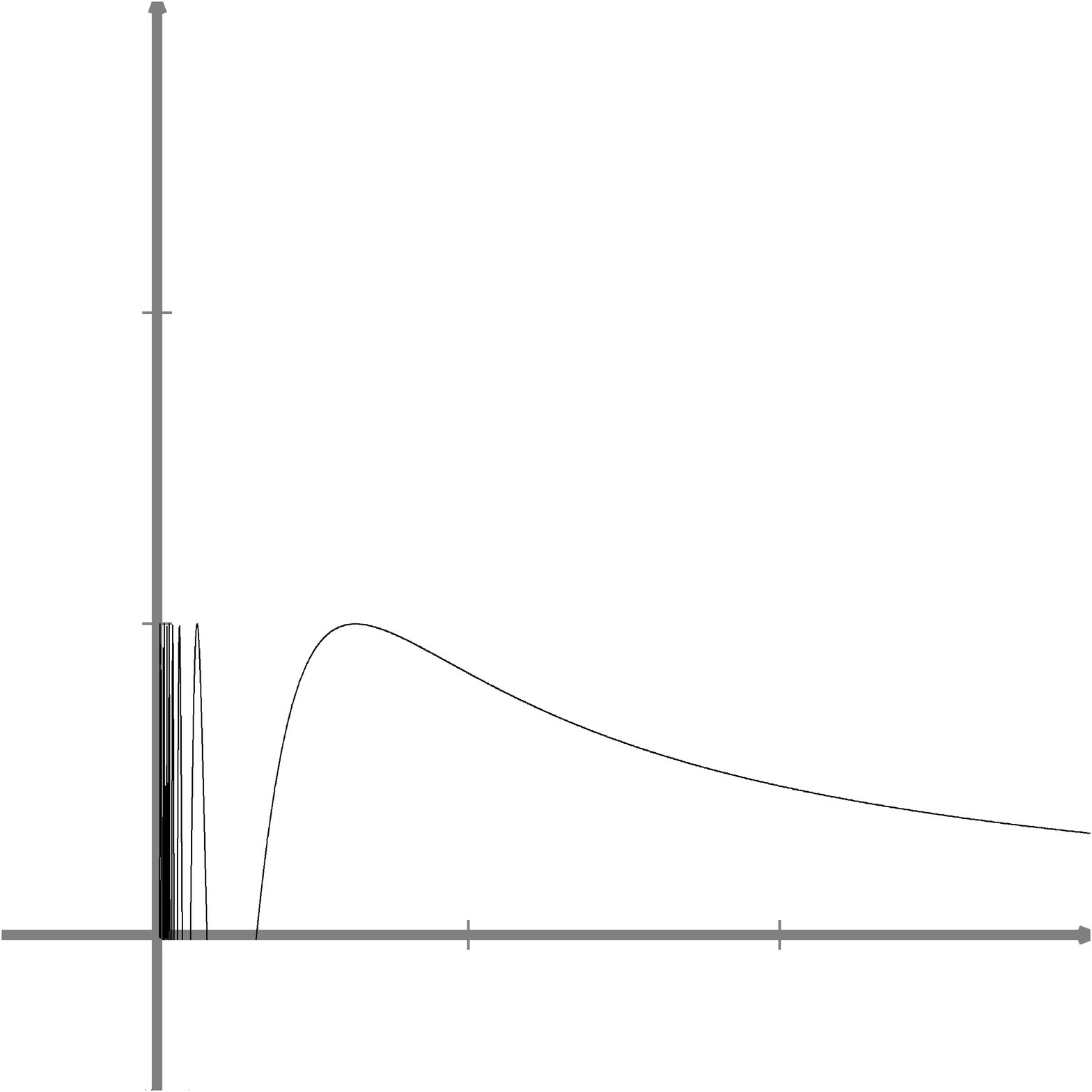}
\hspace{0.5cm}
\includegraphics[width=\figuresize]{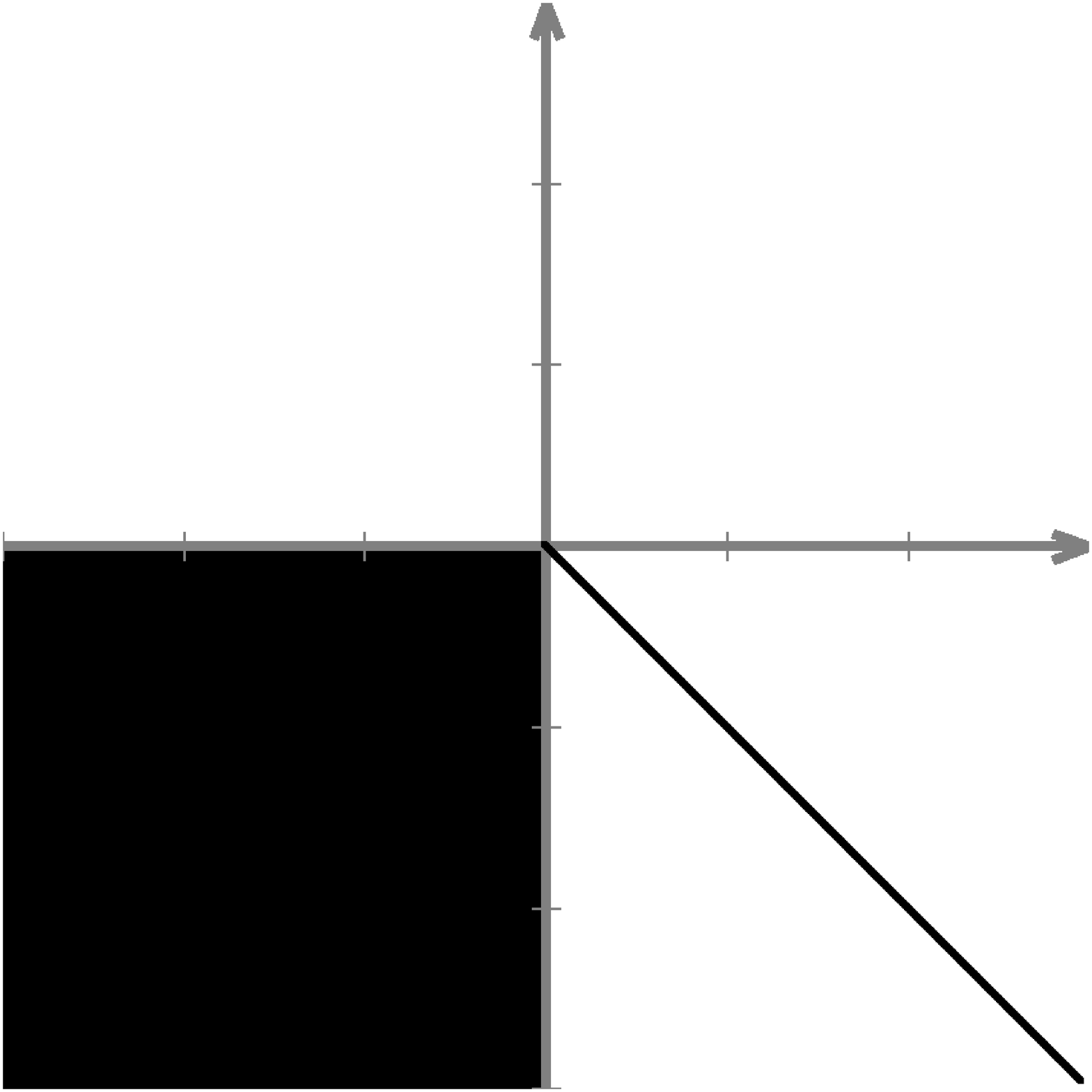}
\caption{$V = \{(x,y) \in {(\erre_{>0})}^2 \ |\  y = \sin \frac{1}{x} \}$ (left picture), then $\ameba_0(V) = \{(x,y) \in \erre^2 \ |\ y\leq 0, x \leq 0 \mbox{ or } x\geq 0, y = -x \}$ (right picture).}  \label{fig:sininv}
\end{center}
\end{figure}

Note that the statement about the dimension can be false for a general set. See figure \ref{fig:sininv} for an example.

\begin{figure}[p]
\begin{center}
\includegraphics[width=\figuresize]{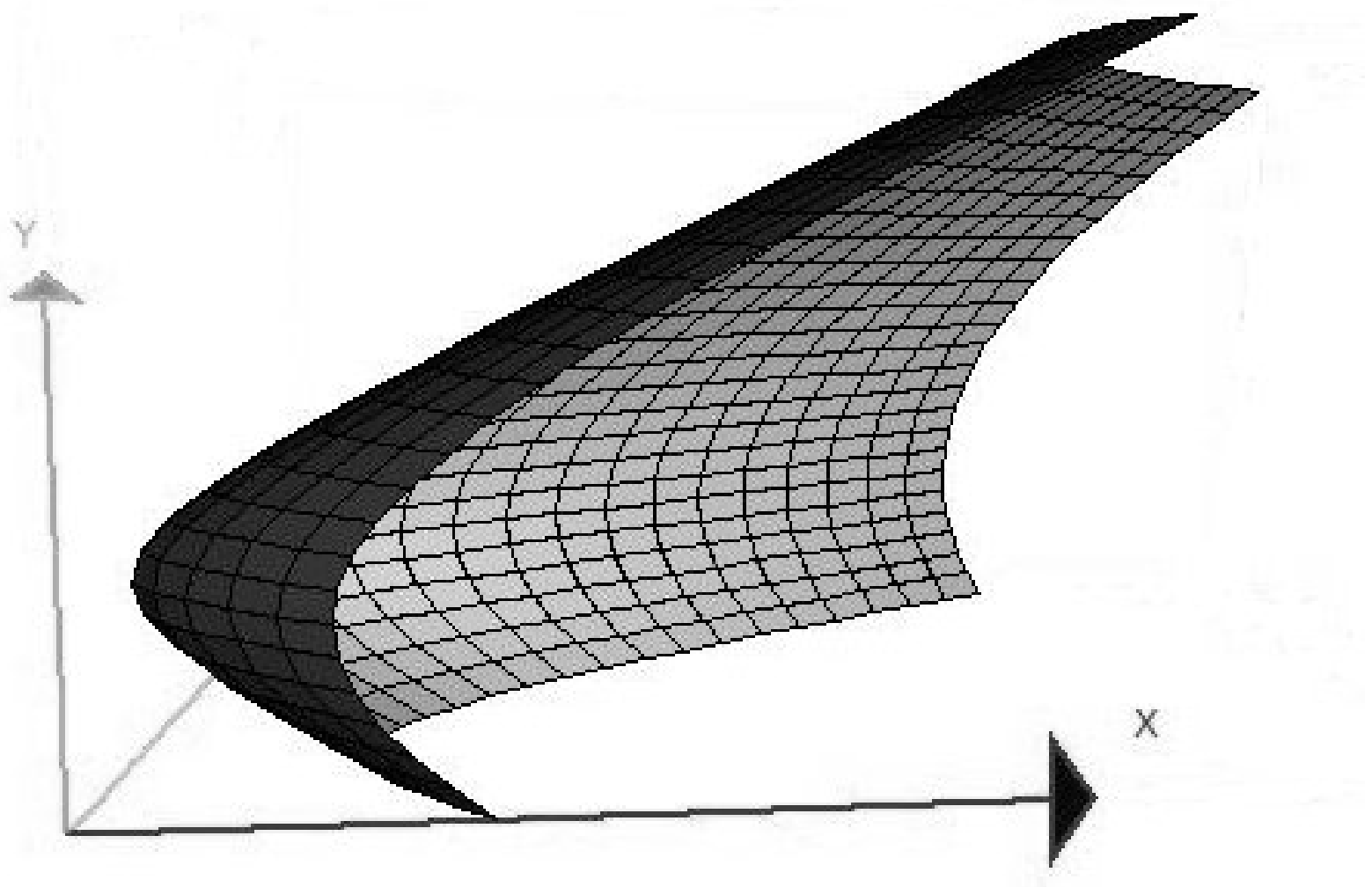}
\hspace{0.5cm}
\includegraphics[width=\figuresize]{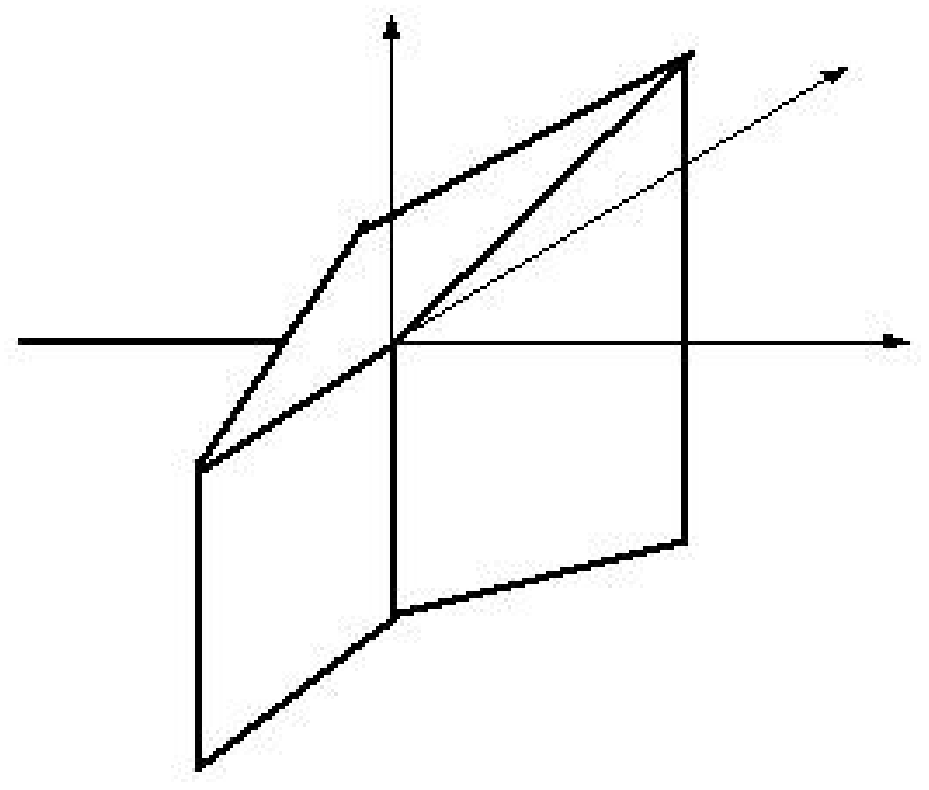}
\caption{$V = \{(x,y,z) \in {(\erre_{>0})}^3 \ |\  (x+1)^2 = 1 + (y-1)^2 + (z-1)^2 = 0 \}$ (left picture), then $\ameba_0(V)$ has an isolated ray along the direction $(-1,0,0)$ and a bi-dimensional part in the half-space $x \geq 0$ (right picture).}  \label{fig:iperboloid}
\end{center}
\end{figure} 

\begin{figure}[p]
\begin{center}
\includegraphics[height=\figuresize]{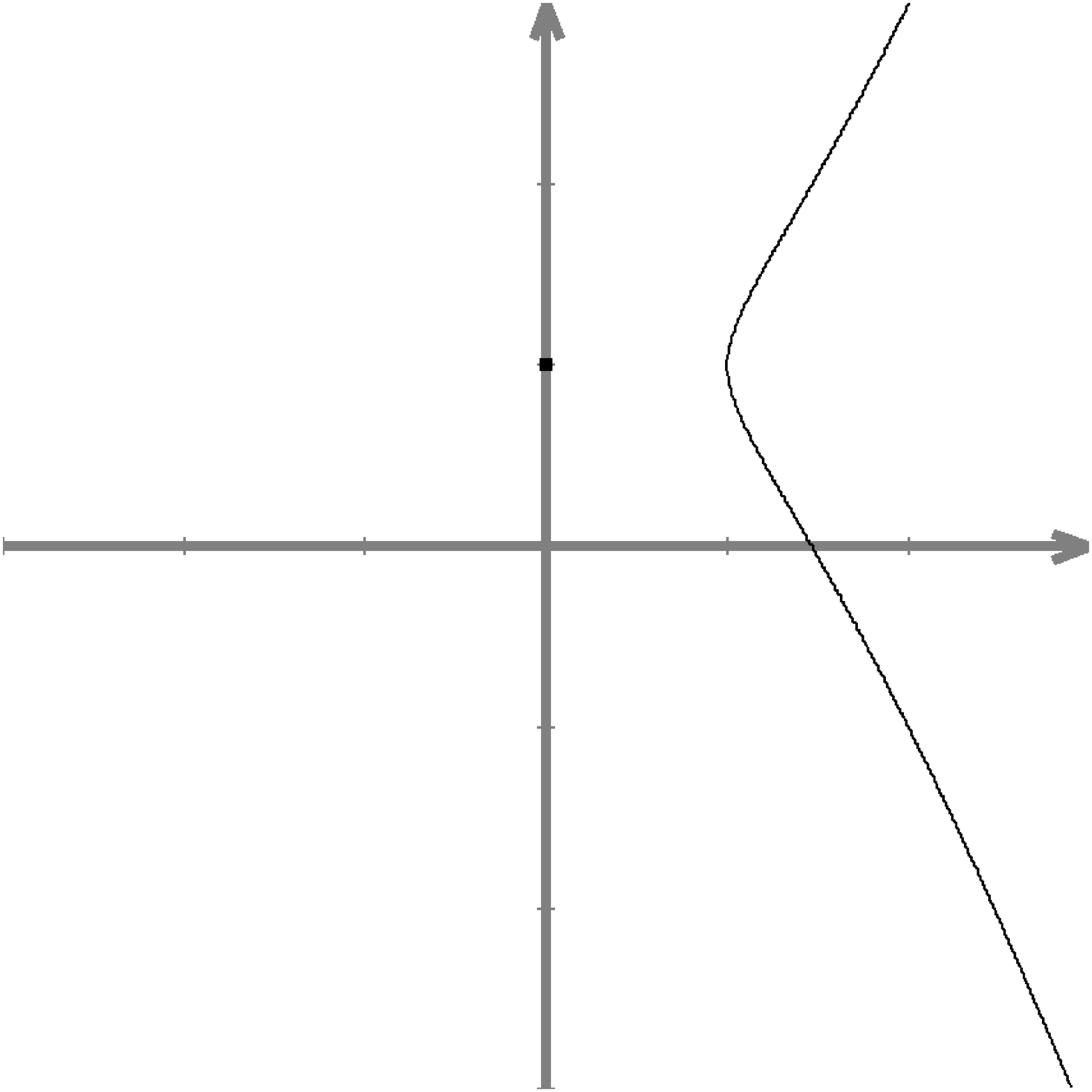}
\hspace{0.5cm}
\includegraphics[width=\figuresize]{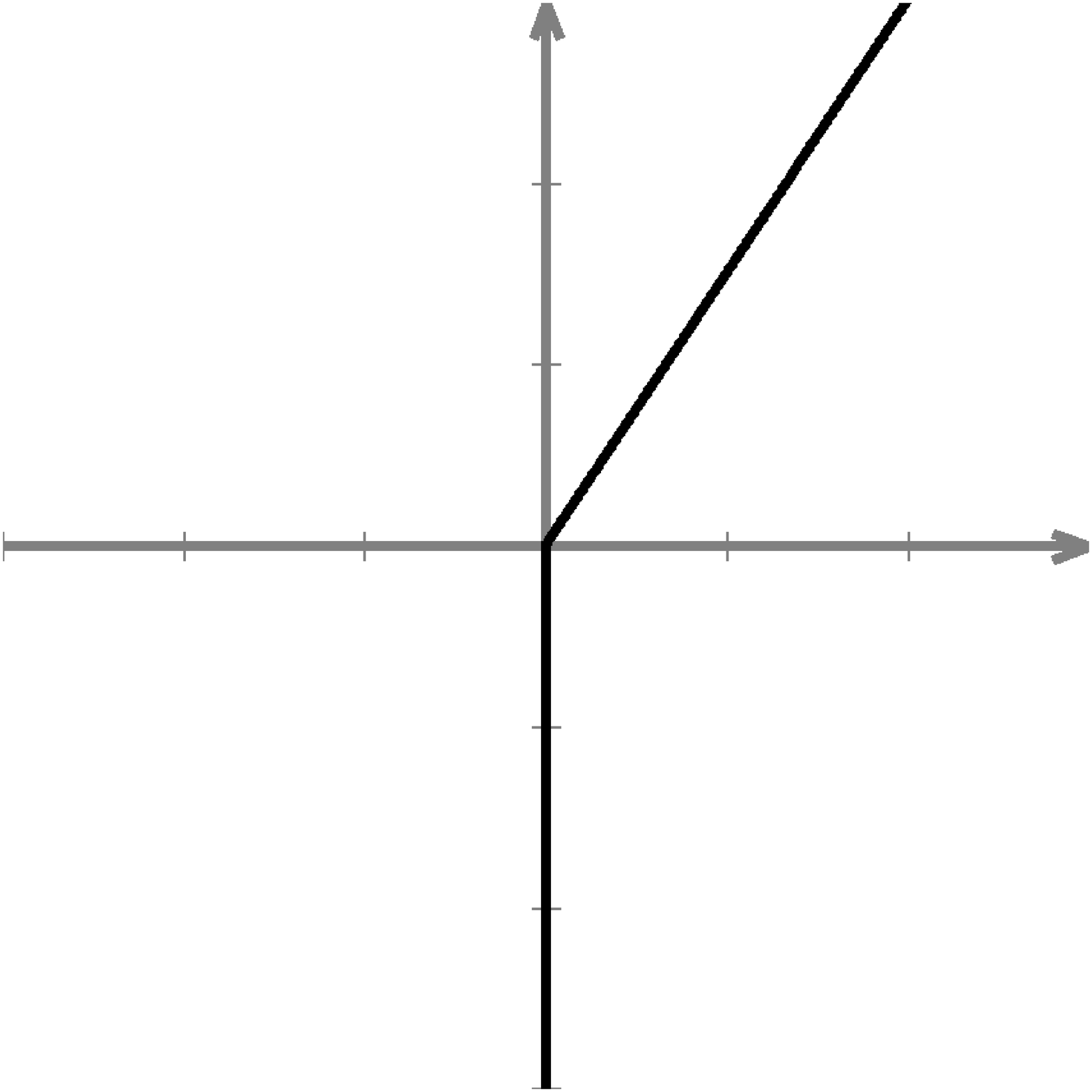}
\caption{$V = \{(x,y) \in \erre^2 \ |\  x^2 + y^2 + 1 = 2y + x^3 \}$ with an isolated point in $(0,1)$ (left picture), then $\ameba_0(V \cap {(\erre_{>0})}^2) = \{(x,y) \in \erre^2 \ |\ x = 0, y \leq 0 \mbox{ or } x\geq 0, 2y = 3x \}$ (right picture).}  \label{fig:cubic}
\end{center}
\end{figure}

Moreover, it is not possible to give more than an inequality, as for every $s \leq m$ it is always possible to find a semi-algebraic set $V \subset {(\erre_{>0})}^m$ such that $\dim V = m$ and $\dim \ameba_0(V) = s$. For example take the parallelepiped $V = [1,2]^{m-s} \times {(\erre_{>0})}^s \subset {(\erre_{>0})}^m$, with $\ameba_0(V) = \{0\}^{m-s} \times {(\erre_{>0})}^s$. It is also possible to find counterexamples of these kind where $V$ is the intersection of ${(\erre_{>0})}^{m+1}$ with an algebraic hypersurface. For example let $S^{m-s} \subset {(\erre_{>0})}^{m-s+1}$ be the sphere with center $(2,\dots,2)$ and radius $1$, then $V = S^{m-s} \times {(\erre_{>0})}^s \subset {(\erre_{>0})}^{m+1}$ has dimension $m$, but $\ameba_0(V) = \{0\}^{m-s+1} \times {(\erre_{>0})}^s$ has dimension $s$.

It is also possible to find a semi-algebraic set $V$ that is the intersection of ${(\erre_{>0})}^n$ with an irreducible pure-dimensional smooth hypersurface, and such that its logarithmic limit set $\ameba_0(V)$ is not pure-dimensional, see for example figure \ref{fig:iperboloid}. Note that the product $V \times S^h$, with $S^h$ the sphere with center $(2,\dots,2)$ and radius $1$ as above, is again the intersection of ${(\erre_{>0})}^{n+h+1}$ with an irreducible pure-dimensional smooth variety, and its logarithmic limit set is lower dimensional and not pure-dimensional.

\section{Non-archimedean description}     \label{sez:non-arch}

\subsection{The Hardy field}   \label{subsez:Hardy}

Let $S$ be a set of symbols expanding $\mathcal{OS}$, and let $\mathfrak{R}=(\erre,a)$ be an o-minimal $S$-structure expanding $\overline{\erre}$ (see subsection \ref{subsez:o-minimal structures} for definitions).

The \nuovo{Hardy field} of $\mathfrak{R}$ can be defined in the following way. If $f,g:\erre_{>0} \freccia \erre$ are two definable functions, we say that they have the same \nuovo{germ} near zero, and we write $f \sim g$, if there exists an  $\varepsilon > 0$ such that $f_{|(0,\varepsilon)} = g_{|(0,\varepsilon)} $. The Hardy field can be defined as the set of germs of definable functions near zero: $H(\mathfrak{R}) = \{f:\erre_{>0} \freccia \erre \ |\ f \mbox{ definable } \} / \sim$. We will denote by $[f]$ the germ of a function $f$.

For every element $a \in \erre$, the constant function with value $a$ defines a germ that is identified with $a$. This defines an an embedding $\erre \freccia H(\mathfrak{R})$. Every relation in the structure $\mathfrak{R}$ defines a corresponding relation on $H(\mathfrak{R})$, and every function in the structure $\mathfrak{R}$ defines a function on $H(\mathfrak{R})$, hence the Hardy field $H(\mathfrak{R})$ can be endowed with an $S$-structure $\overline{H(\mathfrak{R})}$. Given an $(L_S)$-formula $\phi(x_1,\dots,x_n)$, and given definable functions $f_1,\dots,f_n$, we have:
$$\overline{H(\mathfrak{R})}  \vDash \phi([f_1], \dots, [f_n]) \Leftrightarrow \exists \varepsilon >0: \forall t \in (0,\varepsilon): \mathfrak{R} \vDash \phi(f_1(t),\dots,f_n(t))$$

See \cite[sez. 5.3]{Co} for precise definitions and proofs. In particular the $S$-structure $\overline{H(\mathfrak{R})}$ is an elementary extension of the $S$-structure $\mathfrak{R}$. Note that the operations $+$ and $\cdot$ turn $H(\mathfrak{R})$ in a field, the order $\leq$ turn it in a ordered field, and that this field is real closed. Moreover, the $S$-structure $\overline{H(\mathfrak{R})}$ is o-minimal.

Suppose that $S'$ is an expansion of $S$, and that $\mathfrak{R}'$ is an $S'$-structure expanding $\mathfrak{R}$. Then all functions that are definable in $\mathfrak{R}$ are also definable in $\mathfrak{R}'$. This defines an inclusion $H(\mathfrak{R}) \subset H(\mathfrak{R}')$. Note that, by restriction, $\mathfrak{R}'$ has an $S$-structure induced by his $S'$-structure. If $\phi(x_1,\dots,x_n)$ is an $(L_{S})$-formula, and $h_1,\dots,h_n \in H(\mathfrak{R})$, then
$$\overline{H(\mathfrak{R})}  \vDash \phi(h_1, \dots, h_n) \Leftrightarrow \overline{H(\mathfrak{R}')} \vDash \phi(h_1,\dots,h_n)$$
In other words the $S$-structure on $H(\mathfrak{R}')$ is an elementary extension of $\overline{H(\mathfrak{R})}$.

If $\mathfrak{R}$ is polynomially bounded, for every definable function $f$ whose germ is not $0$, there exists $r$ in the field of exponents and $c \in \erre \setminus \{0\}$ such that:
$$\displaystyle \lim_{x \tende 0^+} \dfrac{f(x)}{x^r} = c $$

If $h$ is the germ of $f$, we denote the exponent $r$ by $v(h)$. The map $v:H(\mathfrak{R})\setminus \{0\} \freccia \erre$ is a real valued \nuovo{valuation}, turning $H(\mathfrak{R})$ in a \nuovo{non-archimedean field} of rank one.

The image group of the valuation is the field of exponents of $\mathfrak{R}$, denoted by $\Lambda$. The valuation has a natural section, the map 
$$\Lambda \ni r \freccia x^r \in H(\mathfrak{R})$$

The \nuovo{valuation ring}, denoted by $\ocors$, is the set of all germs bounded in a neighborhood of zero, and the \nuovo{maximal ideal} $m$ of $\ocors$ is the set of all germs infinitesimal in zero. The valuation ring $\ocors$ is \nuovo{convex} with respect to the order $\leq$, hence the valuation topology coincides with the order topology. The map $\ocors \freccia \erre$ sending every element of $\ocors$ in its value in zero, has kernel $m$, hence it identifies in a natural way the residue field $\ocors / m$ with $\erre$.

We will usually denote by $t \in H(\mathfrak{R})$ the germ of the identity function. We have $v(t) = 1$.

As an example we can describe the field $H(\overline{\erre})$. Every element of this field is algebraic over the fraction field $\erre(t)$. Hence $H(\overline{\erre})$ is the real closure of $\erre(t)$, with reference to the unique order such that $t > 0$ and $\forall x \in \erre_{>0}: t < x$. The image of the valuation is $\qu$. Consider the real closed field of formal Puiseaux series with real coefficients, $\erre((t^{\qu})) = \bigcup_{n \geq 1} \erre((t^{1/n}))$. The elements of this field have the form
$$x^r(s(x^{1/n})) $$
where $r \in \ze$ and $s$ is a formal power series. As $\erre(t) \subset \erre((t^{\qu}))$ as an ordered field, then $H(\overline{\erre}) \subset \erre((t^{\qu}))$. The elements of $H(\overline{\erre})$ are the elements of $\erre((t^{\qu}))$ that are algebraic over $\erre(t)$. For these elements the formal power series $s$ is locally convergent.

Another example is the field $H(\erre_{\mathrm an^*})$ (see subsection \ref{subsez:o-minimal structures}). By \cite[thm. B]{DS98}, for every element $h$ of this field, exist a number $r \in \erre$, a formal power series 
$$F = \sum_{\alpha \in \erre_{\geq 0}} c_\alpha X^\alpha$$
and a radius $\delta \in \erre_{>0}$ such that: $c_\alpha \in \erre$, $\{\alpha \ |\ c_\alpha \neq 0\}$ is well ordered, the series $\sum_\alpha |c_\alpha| r^{\alpha} < +\infty$, (hence $F$ is convergent and defines a continuous function on $[0,\delta]$, analytic on $(0,\delta)$) and
$$h = [x^r F(x)]$$

Let $\effe$ be a real closed field extending $\erre$. The convex hull of $\erre$ in $\effe$ is a valuation ring denoted by $\ocors_\leq$. This valuation ring defines a valuation $v:\effe^* \freccia \Lambda$, where $\Lambda$ is an ordered abelian group. We say that $\effe$ is a real closed \nuovo{non-archimedean} field of \nuovo{rank} one extending $\erre$ if $\Lambda$ has rank one as an ordered group, or, equivalently, if $\Lambda$ is isomorphic to an additive subgroup of $\erre$. Hence real closed non-archimedean fields of rank one extending $\erre$ have a real valued valuation (non necessarily surjective) well defined up to a scaling factor. This valuation is well defined when we choose an element $t \in \effe$ with $t>0$ and $v(t)>0$, and we choose a scaling factor such that $v(t) = 1$. Now a valuation $v:\effe \freccia \erre$ is well defined, with image $v(\effe^*) = \Lambda \subset \erre$.

Consider the subfield $\erre(t) \subset \effe$. The order induced by $\effe$ has the property that $t>0$ and $\forall x \in \erre_{>0}: t < x$. Hence $\effe$ contains the real closure of $\erre(t)$ with reference to this order, i.e. $H(\overline{\erre})$. Moreover the valuation $v$ on $\effe$ restricts to the valuation we have defined on $H(\overline{\erre})$, as, if $\ocors_{\leq}$ is the valuation ring of $\effe$, $\ocors_{\leq} \cap H(\overline{\erre})$ is precisely the valuation ring $\ocors$ of $H(\overline{\erre})$. In other words every non-archimedean real closed field of rank one $\effe$ extending $\erre$ is a valued extension of $H(\overline{R})$.

\subsection{Non archimedean amoebas}

Let $\effe$ be a non-archimedean real closed field of rank one extending $\erre$, with a fixed real valued valuation $v:\effe^* \freccia \erre$. By convention, we define $v(0) = \infty$, an element greater than any element of $\erre$. The map
$$\effe \ni h \freccia \parallel h\parallel=\exp(-v(h)) \in \erre_{\geq 0}$$
is a \nuovo{non-archimedean norm}. The \nuovo{component-wise logarithm map} can be defined also on $\effe$, by:
$$\Log:{({\effe}_{>0})}^n \ni (h_1, \dots, h_n) \freccia \left(\log(\parallel h_1\parallel),\dots,\log(\parallel h_n\parallel)\right) \in \erre^n$$
Note that $\log(\parallel h\parallel) = -v(h)$.
If $V \subset {({\effe}_{>0})}^n$, the \nuovo{logarithmic image} of $V$ is the image $\Log(V)$.

Let $S$ be a set of symbols expanding $\mathcal{OS}$, and let $(\effe,a)$ be an $S$-structure expanding the $\mathcal{OS}$-structure on the non-archimedean real closed field of rank one $\effe$ extending $\erre$.
If $V \subset {(\effe_{>0})}^n$ is a definable set in $(\effe,a)$, we call the closure of the logarithmic image of $V$ a \nuovo{non-archimedean amoeba}, and we will write $\ameba(V) = \overline{\Log(V)}$.

The case we are more interested in is when $\mathfrak{R} = (\erre,a)$ is an o-minimal, polynomially bounded $S$-structure expanding $\overline{\erre}$, and $H(\mathfrak{R})$ is the Hardy field, with its natural valuation $v$ and its natural $S$-structure. Non-archimedean amoebas of definable sets of $H(\mathfrak{R})$ are closely related with logarithmic limit sets of definable sets of $\mathfrak{R}$.

Let $\effe \subset \cappa$ be two real closed fields. Let $S$ be a set of symbols expanding $\mathcal{OS}$, let $(\effe,a), (\cappa,b)$ be $S$-structures expanding the $\mathcal{OS}$ structure on the real closed fields and such that $\cappa$ is an elementary extension of $\effe$. Let $V \subset \effe^n$ be a definable set in $(\effe,a)$. We will denote by $\overline{V}$ the \nuovo{extension} of $V$ to the structure $(\cappa,b)$.

For example, if $V \subset \erre^n$ is a definable set in $\mathfrak{R}$, we can always define an extension $\overline{V} \subset {H(\mathfrak{R})}^n$ of $V$ to $H(\mathfrak{R})$.

\begin{lemma} \label{lemma:non-arch amoeba}
Let $\mathfrak{R}$ be a o-minimal polynomially bounded structure. Let $V \subset {(\erre_{>0})}^n$ be a definable set. Then
$$(0,\dots,0,-1) \in \ameba_0(V) \Leftrightarrow (0,\dots,0,-1) \in \Log(\overline{V})$$
\end{lemma}

\begin{proof}
Suppose that $(-1,0,\dots,0) \in \ameba(\overline{V})$. Then there is a point $(x_1, \dots, x_n) \in \overline{V}$ such that $v(x_n) = 1$ and $v(x_i) = 0$ for all $i < n$. Then if $f_1, \dots, f_n$ are definable functions such that $x_i=[f_i]$:
$$\exists \varepsilon > 0 : \forall t \in (0,\varepsilon) : (f_1(t), \dots, f_n(t)) \in V $$
Moreover, when $t \tende 0$ we have that $f_n(t) \tende 0$ and $f_i(t) \tende a_i > 0$ for $i < n$. Hence $V$ contains a sequence tending to $(a_1, \dots, a_{n-1},0)$ with $a_1, \dots a_{n-1} \neq 0$, and $\ameba_0(V)$ contains $(0,\dots,0,-1)$.

Vice versa, suppose that $(0,\dots,0,-1) \in \ameba_0(V)$. Then, by theorem \ref{teo:loglimset} there exists a sequence $x(k)$ in $V$ such that $x(k) \tende (a_1,\dots,a_{n-1},0)$, where $a_1, \dots, a_{n-1} > 0$. Let $\varepsilon$ be a number less than all the numbers $a_1, \dots, a_{n-1}$, and consider the set:
$$\{x \in \erre \ |\ \exists x_1, \dots, x_{n-1} : |x_i-a_i| < {\textstyle\frac{1}{2}}\varepsilon  \mbox{ and } (x_1, \dots, x_{n-1},x) \in V \} $$
As this set is definable, and as it contains a sequence converging to zero, it must contain an interval of the form $(0,\delta)$, with $\delta > 0$. In one formula:
$$\forall x \in (0,\delta): \exists x_1, \dots, x_{n-1} : |x_i-a_i| < {\textstyle\frac{1}{2}}\varepsilon  \mbox{ and } (x_1, \dots, x_{n-1},x) \in V $$
This sentence can be turned into a first order $S$-formula using a definition of $V$. This formula must also hold for $H(\mathfrak{R})$. We can choose an $x \in H(\mathfrak{R})$, with $x > 0$ and $v(x) = 1$. Then $x < \delta$, hence 
$$\exists x_1, \dots, x_{n-1} : |x_i-a_i| < {\textstyle\frac{1}{2}}\varepsilon  \mbox{ and } (x_1, \dots, x_{n-1},x) \in \overline{V} \} $$
Now $v(x_i) = 0$ for all $i > 1$, as $|x_i-a_i| < {\textstyle\frac{1}{2}}\varepsilon$. Hence 
$$\Log(x_1, \dots, x_{n-1},x) = (0,\dots,0,-1)$$
\end{proof}

\begin{theorem} \label{teo:non-arch amoebas 1}
Let $\mathfrak{R}$ be a o-minimal polynomially bounded structure with field of exponents $\Lambda$. Let $V \subset {(\erre_{>0})}^n$ be a definable set. Then 
$$\ameba_0(V) \cap \Lambda^n = \Log(\overline{V}) $$
\end{theorem}

\begin{proof}
We need to prove that for all $x \in \Lambda^n$, $x \in \ameba_0(V) \Leftrightarrow x \in \ameba(\overline{V})$. We choose a matrix $B$ with entries in $\Lambda$ sending $x$ in $(0,\dots,0,-1)$. Then we conclude by the previous lemma applied to the definable set $\overline{B}(V)$.
\end{proof}

\begin{theorem}   \label{teo:density}
Let $\effe \subset \cappa$ be two non-archimedean real closed fields of rank one extending $\erre$, with a choice of a real valued valuation defined by an element $t \in \effe$. Denote the value groups by $\Lambda = v(\effe^*)$ and $\Omega = v(\cappa^*)$. Let $S$ be a set of symbols expanding $\mathcal{OS}$, let $(\effe,a), (\cappa,b)$ be $S$-structures expanding the $\mathcal{OS}$ structure on the real closed fields and such that $\cappa$ is an elementary extension of $\effe$. Let $V$ be a definable set in $(\effe,a)$, and $\overline{V}$ be its extension to $(\cappa,b)$. Then $\Log(V) \subset \Lambda^n$ is dense in $\Log(\overline{V}) \subset \Omega^n$.
\end{theorem}

\begin{proof}
Suppose, by contradiction, that $x \in \Log(\overline{V})$ and it is not in the closure of $\Log(V)$. Then there exists an $\varepsilon > 0$ such that the cube 
$$C = \{ y \in \erre^n \ |\ |y_1-x_1|<\varepsilon, \dots, |y_n-x_n| < \varepsilon  \} $$
does not contain points of $\Log(V)$.

Let $h \in \overline{V}$ be an element such that $\Log(h) = x$, and let $d \in \effe$ be an element such that $0 < v(d) < \varepsilon$.
Consider the cube 
$$E = (\frac{h_1}{d}, h_1 d) \times (\frac{h_2}{d}, h_2 d)\times \dots \times (\frac{h_n}{d}, h_n d) \subset {\cappa}^n $$
The image $\Log(E)$ is contained in $C$, hence $E \cap V$ is empty. But, as $(\cappa,b)$ is an elementary extension of $(\effe,a)$, also $E \cap \overline{V}$ is empty. This is a contradiction as $h \in E$ and $h \in \overline{V}$.
\end{proof}

\begin{corollary}
Let $S$ be a set of symbols expanding $\mathcal{OS}$, and let $\mathfrak{R}=(\erre,a)$ be an o-minimal polynomially bounded $S$-structure with field of exponents $\Lambda$, expanding $\overline{\erre}$. Let $V \subset {(\erre_{>0})}^n$ be a definable set. Suppose that there exists a subfield $\Omega \subset \erre$ such that $\Lambda \subset \Omega$ and $\mathfrak{R}^\Lambda$ is o-minimal and polynomially bounded. Then $\ameba_0(V) \cap \Lambda^n$ is dense in $\ameba_0(V) \cap \Omega^n$.
\end{corollary}

\begin{proof}
Consider the Hardy fields $H(\mathfrak{R})$ and $H(\mathfrak{R}^\Lambda)$. We denote by $\overline{V}$ the extension of $V$ to $H(\mathfrak{R})$, and by $\overline{\overline{V}}$ the extension of $V$ to $H(\mathfrak{R}^\Lambda)$. By theorem \ref{teo:non-arch amoebas 1} $\ameba_0(V) \cap \Lambda^n = \Log(\overline{V})$ and $\ameba_0(V) \cap \Omega^n = \Log(\overline{\overline{V}})$. As we said above, the $S$-structures on $H(\mathfrak{R})$ and $H(\mathfrak{R}^\Lambda)$ are elementary equivalent. The statement follows by the previous theorem.
\end{proof}

\begin{corollary}   \label{corol:log e nonarch}
Let $S$ be a set of symbols expanding $\mathcal{OS}$, and let $\mathfrak{R}=(\erre,a)$ be a regular polynomially bounded $S$-structure with field of exponents $\Lambda$. Let $V \subset {(\erre_{>0})}^n$ be a set that is definable in $\mathfrak{R}$. We denote by $\overline{V}$ the extension of $V$ to $H(\mathfrak{R})$ and by $\overline{\overline{V}}$ the extension of $V$ to $H(\mathfrak{R}^\erre)$. Then
$$\ameba_0(V) = \Log(\overline{\overline{V}})$$ 
Moreover the subset $\ameba_0(V) \cap \Lambda^n$ is dense in $\ameba_0(V)$, and, as $\ameba_0(V)$ is closed,
$$\ameba(\overline{V}) = \ameba(\overline{\overline{V}}) = \Log(\overline{\overline{V}})$$
\end{corollary}

\begin{corollary}    \label{corol:non-arch semi-alg}
Let $V \subset {(\erre_{>0})}^n$ be a semi-algebraic set. Then $\ameba_0(V) \cap \qu^n$ is dense in $\ameba_0(V)$. Let $\effe$ be a non-archimedean real closed field of rank one extending $\erre$, and let $\overline{V}$ be the extension of $V$ to $\effe$. Then
$$\ameba_0(V) = \ameba(\overline{V}) $$
If $\effe$ extends $H({\overline{\erre}}^{\erre})$, then
$$\ameba(\overline{V}) = \Log(\overline{V})$$
\end{corollary}

As a further corollary, we prove the following proposition, that will be needed later.

\begin{proposition}    \label{prop:projection}
Let $V \subset ({\erre_{>0}})^n$ be a set definable in a regular polynomially bounded structure, and let $\pi:\erre^n \freccia \erre^m$ be the projection on the first $m$ coordinates (with $m < n$). Then we have
$$\ameba_0(\pi(V)) = \pi(\ameba_0(V)) $$
\end{proposition}

\begin{proof}
It follows easily from corollary \ref{corol:log e nonarch} and from the fact that $\overline{\pi}(\overline{V}) = \overline{\pi(V)}$.
\end{proof}

\subsection{Patchworking families}  \label{subsez:patchworking}

Let $S$ be a set of symbols expanding $\mathcal{OS}$, and let $\mathfrak{R} = (\erre,a)$ be an $S$-structure expanding $\overline{\erre}$. If $V \subset {(H(\mathfrak{R})_{>0})}^n$ is definable, there exists a first order $S$-formula $\phi(x_1, \dots, x_n, y_1, \dots, y_m)$, and parameters $a_1, \dots, a_m \in H(\mathfrak{R})$ such that
$$V = \{(x_1, \dots, x_n) \ |\ \phi(x_1, \dots, x_n, a_1, \dots, a_m)\} $$
Choose definable functions $f_1, \dots, f_m$ such that $[f_i] = a_i$. These data defines a definable set in $\mathfrak{R}$:
$$\widetilde{V} = \{(x_1, \dots, x_n, t) \in {(\erre_{>0})}^{n+1} \ |\ \phi(x_1, \dots, x_n, f_1(t), \dots, f_m(t)) \}$$
Suppose that $\phi'(x_1, \dots, x_n, y_1, \dots, y_{m'})$ is another formula defining $V$ with parameters $a_1', \dots, a_{m'}'$, and that $f_1', \dots, f_{m'}'$ are definable functions such that $[f_i']=a_i$. These data defines:
$$\widetilde{V}' = \{(x_1, \dots, x_n, t) \in {(\erre_{>0})}^{n+1} \ |\ \phi'(x_1, \dots, x_n, f_1'(t), \dots, f_{m'}'(t)) \}$$
As both formulae defines $V$ we have:
$$\overline{H(\mathfrak{R})} \vDash \forall x_1, \dots, x_n: \phi(x_1, \dots, x_n, a_1, \dots, a_m) \Leftrightarrow \phi'(x_1, \dots, x_n, a_1', \dots, a_{m'}')$$
As we said above, we have:
$$\exists \varepsilon >0: \forall t \in (0,\varepsilon): \mathfrak{R} \vDash \forall x_1, \dots, x_n: $$
$$\phi(x_1, \dots, x_n, f_1(t), \dots, f_m(t)) \Leftrightarrow \phi'(x_1, \dots, x_n, f_1'(t), \dots, f_{m'}'(t))$$
Hence 
$$\widetilde{V} \cap \left(\erre^n \times (0,\varepsilon)\right) = \widetilde{V}' \cap \left(\erre^n \times (0,\varepsilon)\right)$$
and the set $\widetilde{V}$ is ``well defined for small enough values of $t$''.
Actually we prefer to see the set $\widetilde{V}$ as a parametrized family:
$$V_t = \{ (x_1, \dots, x_n)\in {(\erre_{>0})}^n \ |\ (x_1, \dots, x_n,t) \in \widetilde{V} \}$$
we can say that the set $V$ determines the germ near zero of this parametrized family. We will use the notation $V_* = {(V_t)}_{t>0}$ for the family, and we will call these families \nuovo{patchworking families} determined by $V$, as they are a generalization of the patchworking families of $\cite{Vi}$.

Given a patchworking family $V_*$, we can define the \nuovo{tropical limit} of the family as:
$$\ameba_0(V_*) = \lim_{t \freccia 0} \ameba_t(V_t) = \lim_{t \freccia 0} \Log_{\frac{1}{t}}(V_t)$$
This is a closed subset of $\erre^n$. Note that this set only depends on $V$. If $V$ is the extension to $H(\mathfrak{R})$ of a definable subset $W \subset \erre^n$, then the patchworking family $V_t$ is constant: $V_t = W$, and the tropical limit is simply the logarithmic limit set: $\ameba_0(V_*) = \ameba_0(W)$.

Consider the logarithmic limit set of $\widetilde{V}$:
$$\ameba_0(\widetilde{V}) = \lim_{t \freccia 0} \ameba_t(\widetilde{V}) \subset \erre^{n+1}$$
As in subsection \ref{subsez:polyh struc}, we consider the set
$$H = \{(x_1, \dots, x_n) \in \erre^n \ |\ (x_1, \dots, x_n,-1) \in \ameba_0(\widetilde{V}) \}$$

Note that 
$$\Log_{\left(\frac{1}{t}\right)}^{-1}(\erre^n \times \{-1\}) = {(\erre_{>0})}^n \times \{t\} $$
Hence $\ameba_0(V_*) = \lim_{t \freccia 0} \Log_{\frac{1}{t}}(V_t) = H$.

Now consider the extension of the set $\widetilde{V}$ to the Hardy field $H(\mathfrak{R})$, we denote it by $\overline{\widetilde{V}}$. By the results of the previous section, we know that $\ameba(\overline{\widetilde{V}}) = \ameba_0(\widetilde{V})$. If we denote by $t\in H(\mathfrak{R})$ the germ of the identity function, we have that
$$V = \{(x_1, \dots,x_n) \ |\ (x_1, \dots, x_n,t) \in \overline{\widetilde{V}}\} $$
as, for $i \in \{1, \dots, m\}$, we have $f_i(t) = a_i$. Hence, as $\log|t| = -1$, $\ameba(V) \subset H = \ameba_0(V_*)$.

\begin{lemma}
$(0,\dots,0) \in \ameba_0(V_*) \Leftrightarrow (0,\dots,0) \in \Log(V)$.
\end{lemma}

\begin{proof}
$\Rightarrow$: This follows from what we said above.

$\Leftarrow$: It follows from the second part of the proof of lemma \ref{lemma:non-arch amoeba}, applied to the set $\widetilde{V}$.
\end{proof}

Let $\lambda \in \Lambda^n$. We define a twisted set 
$$V^\lambda = \{x \in {H(\mathfrak{R})}^n \ |\ \phi(t^{-\lambda_1}x_1, \dots, t^{-\lambda_n}x_n, a_1, \dots, a_m)\}$$
Then $\lambda \in \Log(V) \Leftrightarrow (0,\dots,0) \in \Log(V^{\lambda})$. Then we define
$$H^{\lambda} = \{(x_1, \dots, x_n) \in \erre^{n} \ |\ (x_1, \dots, x_n,-1) \in \ameba_0(\widetilde{V^{\lambda}}) \}$$
Now $H^{\lambda}$ is simply $H$ translated by the vector $-\lambda$. Hence we get:

\begin{lemma} 
For all $\lambda \in \Lambda$, we have $\lambda \in \ameba_0(V_*) \Leftrightarrow \lambda \in \Log(V)$.
\end{lemma}

Using these facts we can extend the results of the previous sections about logarithmic limit sets and their relations with non-archimedean amoebas, to tropical limits of patchworking families. For example we can prove the following statements.

\begin{theorem}
Let $S$ be a structure expanding $\mathcal{OS}$, and let $\mathfrak{R} = (\erre,a)$ be a regular polynomially bounded $S$-structure with field of exponents $\Lambda$. Let $V$ be a definable subset of the Hardy field $H(\mathfrak{R})$, and let $V_*$ be a patchworking family determined by $V$. Then the following facts hold:
\begin{enumerate}
\item $\ameba_0(V_*)$ is a polyhedral complex with dimension less than or equal to the dimension of $V$.
\item $\ameba_0(V_*) \cap \Lambda^n = \Log(V)$.
\item $\ameba_0(V_*) = \ameba(V)$.
\item $\ameba_0(V_*) \cap \Lambda^n$ is dense in $\ameba_0(V_*)$.
\end{enumerate}
\end{theorem}

\begin{proof}
Every statement follows from the corresponding statement about logarithmic limit sets, and from the facts exposed above.
\end{proof}

For every point $\lambda \in \Lambda$, the twisted set $V^{\lambda}$ defines a germ of patchworking family $V^{\lambda}_*$. The limit
$$V^{\lambda}_0 = \lim_{t \tende 0} V^{\lambda}_t$$
is a definable set in $\erre^n$ and it has the properties of the set $W$ of subsection \ref{subsez:polyh struc}. The difference is that now the set $V^{\lambda}_0$ is well defined, and it depends only on $\lambda$.

\begin{theorem}
Let $S$ be a set of symbols expanding $\mathcal{OS}$, and let $\mathfrak{R} = (\erre,a)$ be an $S$-structure expanding $\overline{R}$, that is o-minimal and polynomially bounded, with field of exponents $\Lambda$. Let $V$ be a definable subset of the Hardy field $H(\mathfrak{R})$. Then we have
$$\forall \lambda \in \Lambda^n: \lambda \in \ameba(V) \Leftrightarrow V^{\lambda}_0  \neq \emptyset$$
Moreover, if $\Lambda = \erre$, for all $\lambda \in \erre$, there exists a neighborhood $U$ of $\lambda$ in $\ameba(V)$ such that the translation of $U$ by $-\lambda$ is a neighborhood of $(0,\dots,0)$ in $\ameba_0(V^{\lambda}_0)$.
\end{theorem}

\begin{proof}
It follows from the arguments above and from theorem \ref{teo:neighborhood}.
\end{proof}

The set $V^{\lambda}_0$ can be called \nuovo{initial set} of $\lambda$, as it plays the role of the initial ideal of \cite{SS04}. The difference is that $V^{\lambda}_0$ is a geometric object, while the initial ideal of \cite{SS04} is a combinatorial one.

\section{Comparison with other constructions} \label{sez:comp}

\subsection{Complex algebraic sets}   \label{subsez:complex}

Logarithmic limit sets of complex algebraic sets are a particular case of logarithmic limit sets of real semi-algebraic sets, in the following sense. Let $V \subset \ci^n$ be a complex algebraic set, and consider the real semi-algebraic set
$$|V| = \{ x \in {(\erre_{>0})}^n \ |\ \exists z \in V: |z|=x \} $$

The logarithmic limit set of $V$ as defined in \cite{Be71} is precisely the logarithmic limit set of $|V|$ in our notation. Hence all the results we got about logarithmic limit sets of real semi-algebraic sets produce an alternative proof of the same results for complex algebraic sets, that were originally proved partly in \cite{Be71} and partly in \cite{BG81}.

Even the description of logarithmic limit sets via non-archimedean amoebas can be translated to complex algebraic sets. Let $\effe$ be a non-archimedean real closed field of rank one extending $\erre$, and let $v$ be a choice of a real valued valuation on $\effe$, as in subsection \ref{subsez:Hardy}. The field $\cappa = \effe[i]$ is an algebraically closed field extending $\ci$, with an extended valuation $v:\cappa^* \freccia \erre$ defined by $v(a+bi) = \min(v(a), v(b))$. The \nuovo{component-wise logarithm map} can be extended to $\cappa$, by:
$$\Log:{\cappa}^n \ni (z_1, \dots, z_n) \freccia (-v(z_1),\dots,-v(z_n)) \in \erre^n$$

On $\cappa$ there is also the \nuovo{complex norm} $|\cdot|:\cappa \freccia \effe_{\geq 0}$ defined by $|a+bi| = \sqrt{a^2 + b^2}$. Now if $V$ is an algebraic set in $\cappa^n$, the set
$$|V| = \{ x \in {(\effe_{>0})}^n \ |\ \exists z \in V: |z|=x \} $$
is a semi-algebraic set in $\effe^n$. The \nuovo{logarithmic image} of $V$ is the image $\Log(V)$, and the non-archimedean amoeba $\ameba(V)$ is the closure of this image. As expected, $\Log(V)=\Log(|V|)$ and $\ameba(V)=\ameba(|V|)$. Moreover, if $V \subset \ci^n$ is an algebraic set, and $\overline{V} \subset \cappa^n$ is its extension to $\cappa$, then $\overline{|V|} = |\overline{V}|$. These facts directly give the relation between logarithmic limit sets of complex algebraic sets and non-archimedean amoebas in algebraically closed fields.

The same relation holds with patchworking families. Let $\mathfrak{R} = (\erre,a)$ be a regular polynomially bounded structure with field of exponents $\Lambda$, let $\effe = H(\mathfrak{R})$ and $\cappa = H(\mathfrak{R})[i]$ and let $V \subset \cappa^n$ be an algebraic set. There are polynomials $f_1, \dots, f_m \in \cappa[x_1, \dots, x_n]$ such that $V = V(f_1, \dots, f_m)$. Every polynomial $f_j$ has the form
$$f_j = \sum_{\omega \in \ze^n} (a_{j,\omega} + i b_{j,\omega}) x^\omega$$
where $a_{j,\omega},b_{j,\omega} \in H(\mathfrak{R})$. Choose representatives functions $\alpha_{j,\omega},\beta_{j,\omega}$ such that $[\alpha_{j,\omega}] = a_{j,\omega},[\beta_{j,\omega}] = b_{j,\omega}$. This choice defines families of polynomials
$$f_{j,t} = \sum_{\omega \in S(f)} (\alpha_{j,\omega}(t) + i \beta_{j,\omega}(t)) x^\omega$$
and a corresponding family of algebraic sets in $\ci^n$
$$V_t = V(f_{1,t}, \dots, f_{m,t})$$
We will call these families \nuovo{patchworking families} because they generalize the patchworking polynomial of \cite[Part 2]{Mi}, and we will denote the family by $V_* = (V_t)$. The family $V_*$ depends of the choice of the polynomials $f_j$ and of the definable functions $\alpha_{j,\omega}, \beta_{j,\omega}$. If we change these choices we get another patchworking family coinciding with $V_*$ for $t \in (0,\varepsilon)$. The \nuovo{tropical limit} of one such family is
$$\ameba_0(V_*) = \lim_{t \freccia 0} \ameba_t(V_t)$$
As before, $|V|$ is a semi-algebraic set in ${H(\mathfrak{R})}^n$, and if $|V|_* = (|V|_t)$ is a patchworking family defined by $|V|$, then there exists an $\varepsilon > 0$ such that for $t \in (0,\varepsilon)$ we have $|V_t| = |V|_t$. Hence we have that $$\ameba_0(V_*) = \ameba_0(|V|_*)$$
and we can get the properties of the tropical limit of complex patchworking families as a corollary of the properties of tropical limits of real patchworking families.

\begin{figure}[htbp]
\begin{center}
\includegraphics[height=\figuresize]{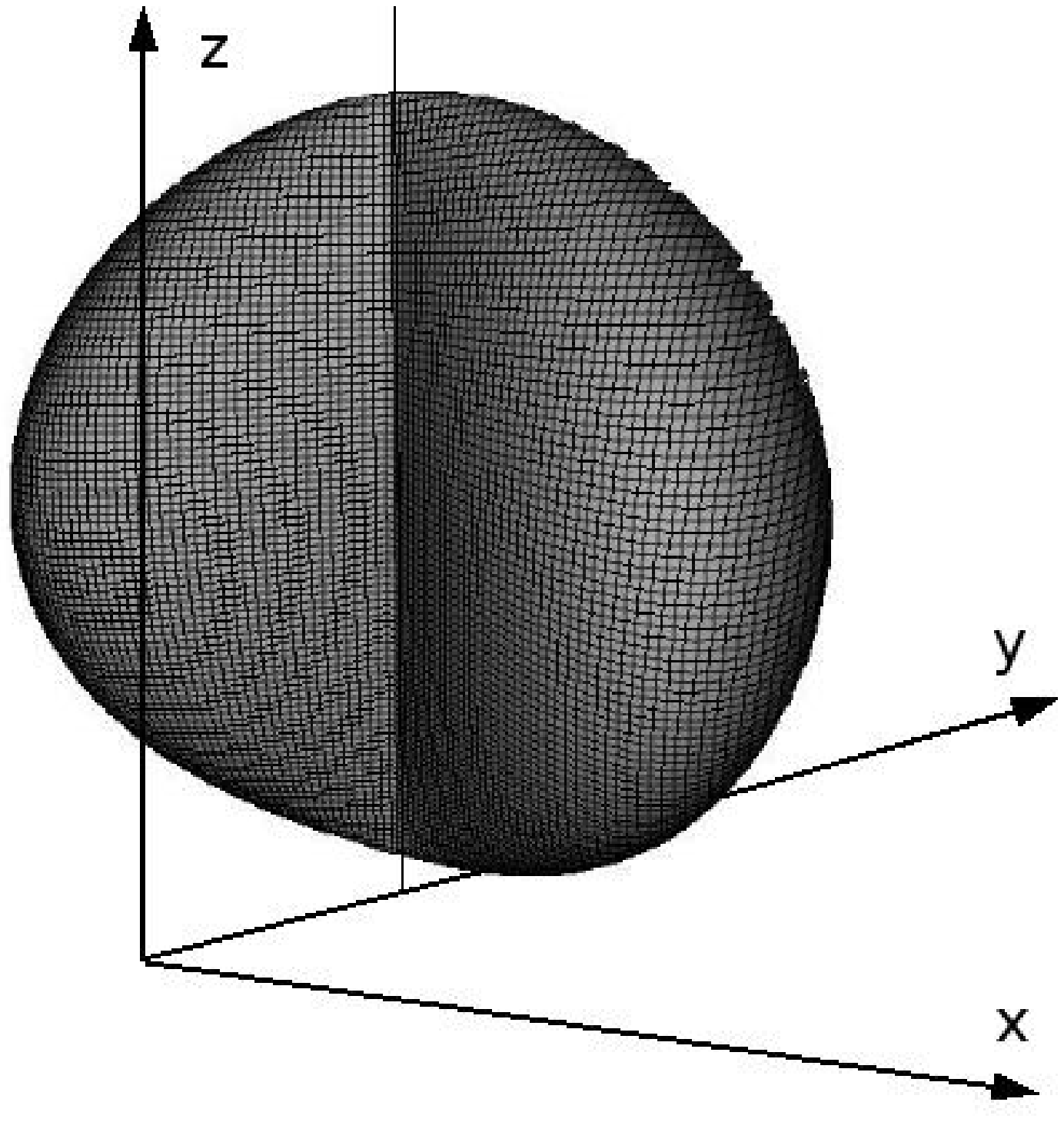}
\hspace{0.5cm}
\includegraphics[width=\figuresize]{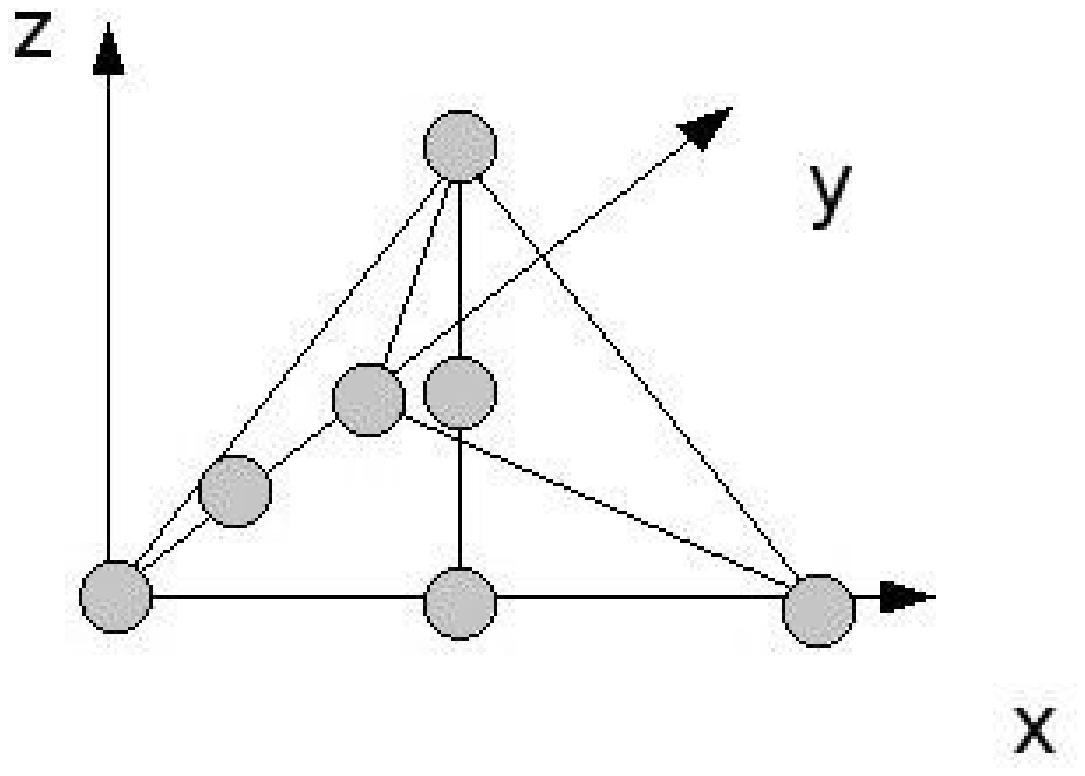}
\caption{$V = \{(x,y,z) \in \erre^3 \ |\  x^2(1-(z-2)^2) = x^4 + (y-1)^2 \}$, it is an irreducible surface, but it has a ``stick'', the line $\{y=1,x=0\}$. The logarithmic limit set of $V \cap {(\erre_{>0})}^3$ is only the ray in the direction $(-1,0,0)$, but this ray is contained in the interior part of a face of the dual fan of the newton polytope of the defining polynomial $x^2(1-(z-2)^2) - x^4 - (y-1)^2$. }  \label{fig:umbrella}
\end{center}
\end{figure}

Let $f \in \erre[x_1, \dots, x_n]$. Let $V$ be the intersection of the zero locus of $f$ and ${(\erre_{>0})}^n$, and let $V_\ci$ be the zero locus of $f$ in $\ci^n$. As $V \subset V_\ci$, the logarithmic limit set of $V$ is included in the logarithmic limit set of $V_\ci$. Moreover, as $V_\ci$ is a complex hypersurface, it is possible to give an easy combinatorial description of $\ameba_0(V_\ci)$, it is simply the dual fan of the newton polytope of $f$. Unfortunately, it is not possible, in general, to use this fact to understand the combinatorics of $\ameba_0(V)$. There are examples where $V$ is an irreducible hypersurface, and $\ameba_0(V)$ is a subpolyhedron of $\ameba_0(V_\ci)$ that is not a subcomplex. For example, if $f$ is as in figure \ref{fig:umbrella}, the logarithmic limit set of $V$ is only the ray in the direction $(-1,0,0)$, but this ray lies in the interior of a face of $\ameba_0(V_\ci)$.

\subsection{Positive tropical varieties}

In this subsection we compare the notion of non-archimedean amoebas for real closed fields that we studied in this chapter with a similar object called positive tropical variety studied in \cite{SW}.  

To be consistent with \cite{SW}, we will denote by $\cappa = \bigcup_{n=1}^\infty \ci((t^{1/n}))$ the algebraically closed field of formal Puiseaux series with complex coefficients, whose set of exponents is an arithmetic progression of rational numbers, and by $\effe = \bigcup_{n=1}^\infty \erre((t^{1/n}))$ the subfield of series with real coefficients. $\cappa$ is the algebraic closure of $\effe$. These fields have a natural valuation $v:\cappa \freccia \qu$, with valuation ring $\ocors$, and residue map $r:\ocors \freccia \ci$. Note that the valuation $v$ is compatible with the order of $\effe$, i.e. the valuation ring $\ocors \cap \effe$ is convex for the order, and that $r(\ocors \cap \effe) = \erre$.

We will denote by $\effe_{>0}$ the set of positive elements of the field $\effe$. Following \cite{SW} we will also use the notation:

$$\effe_+ = \{ z \in \cappa \ |\ r\left(\frac{z}{t^{v(z)}}\right) \in \erre_{>0} \} $$

Let $V$ be an algebraic set in $\cappa^n$. The set
$$V_{>0} = V \cap {(\effe_{>0})}^n $$
is a semi-algebraic set, whose non-archimedean amoeba $\ameba(V_{>0})$ (i.e. the closure of the logarithmic image $\Log(V_{>0})$) has been studied in subsection \ref{subsez:patchworking}. In \cite{SW} a similar definition is given. The positive part of $V$ is
$$V_+ = V \cap {(\effe_+)}^n $$
The closure of $\Log(V_+)$ is called \nuovo{positive tropical variety}, and it is denoted by $\mbox{Trop}^+(V)$. From the definition it is clear that $\ameba(V_{>0}) \subset \mbox{Trop}^+(V)$. 

In many examples the sets $\ameba(V_{>0})$ and $\mbox{Trop}^+(V)$ coincide, but it is also possible to construct examples where the inclusion is strict. For example
$$V = \{ z \in \cappa^2 \ |\ x_1^2 + {(x_2-1)}^2 - x_1^3\} $$
Then $V_{>0}$ is the extension to $\effe$ of the set in figure \ref{fig:cubic}.
$$\ameba(V_{>0}) = \{(x_1,x_2) \in \erre^2 \ |\ x_1 = 0, x_2 \leq 0 \mbox{ or } x_1\geq 0, 2x_2 = 3x_1 \}$$
$$\mbox{Trop}^+(V) = \ameba(V_{>0}) \cup \{x_2 = 0, x_1 \leq 0\}$$

A more interesting example where $\ameba(V_{>0}) \subsetneq \mbox{Trop}^+(V)$ is the set 
$$V = \{(x,y,z) \in \cappa^3 \ |\  x^2(1-(z-2)^2) = x^4 + (y-1)^2 \}$$
Here $V_{>0}$ is the extension to $\effe$ of the set in figure \ref{fig:umbrella}. Now $\ameba(V_{>0})$ is just the ray in the direction $(-1,0,0)$. This ray is in the interior part of a face of $\mbox{Trop}^+(V)$. Hence not only the two sets does not coincide, but $\ameba(V_{>0})$ is not a polyhedral subcomplex of $\mbox{Trop}^+(V)$.

\section{Tropical description}           \label{sez:trop}

\subsection{Maslov dequantization}

Every real number $t \in (0,1)$ defines an analytic function:

$$\erre_{>0}\ni z \freccia \log_{\left(\frac{1}{t}\right)} z = \left(\frac{-1}{\log t}\right) \log z \in \erre$$

This function is bijective, with inverse $x \freccia t^{-x}$, and it preserves the order $\leq$. The operations (`$+$' and `$\cdot$') are transformed via conjugation in the following way:

$$x \oplus_t y = \log_{\left(\frac{1}{t}\right)}(t^{-x} + t^{-y})$$ 

$$x \odot_t y = \log_{\left(\frac{1}{t}\right)}(t^{-x} \cdot t^{-y}) = x + y$$

Hence every $t$ induces an $\mathcal{OS}$-structure on $\erre$: 
$$\erre^t = (\erre, \{\leq\},\{\oplus_t,\odot_t\}, \emptyset)$$ 
This structure is isomorphic to $\erre_{>0}$, hence it is an ordered semifield.

In the limit for $t$ tending to zero we have:

$$ \lim_{t \tende 0^+} x \oplus_t y = \max(x,y)$$

The limit $\mathcal{OS}$-structure is called the \nuovo{tropical semifield}:
$$\erre^{trop} = (\erre, \{\leq\},\{\max,+\}, \emptyset)$$ 
This is again an ordered semifield, we will denote its operations by $\oplus = \max$ and $\odot = +$. Note the inequality
$$ x_1 \oplus \dots \oplus x_n \leq x_1 \oplus_t \dots \oplus_t x_n \leq x_1 \oplus \dots \oplus x_n + \log_{\left(\frac{1}{t}\right)} n$$
In other words the convergence of the family $\erre^t$ to the structure $\erre^{trop}$ is uniform. This construction is usually called \nuovo{Maslov dequantization}.

Note that if $\alpha \in \erre_{>0}$, the function
$$\erre_{>0} \ni x \freccia x^\alpha \in \erre_{>0}$$
is transformed, via conjugation with the map $\log_{\left(\frac{1}{t}\right)}$, in the map:
$$\erre \ni x \freccia \log_{\left(\frac{1}{t}\right)}\left({(t^{-x})}^\alpha\right) = \alpha x$$
As this map does not depend on $t$, it induces also a map in the limit structure $\erre^{trop}$. With these maps, $\erre^t$ and $\erre^{trop}$ become $\mathcal{OS}^\erre$-structures.

The family of maps $\Log_t$, which we used to construct the logarithmic limit sets, is the Maslov dequantization applied coordinate-wise to ${(\erre_{>0})}^n$.

\subsection{Dequantization of formulae}

An $L_{\mathcal{OS}^\erre}$-term $u(x_1, \dots, x_n, y_1, \dots y_m)$ (see \cite[Chap. II, def. 3.1]{EFT84})   and the constants $a_1, \dots, a_m \in \erre_{>0}$ define a function:
$$U: {(\erre_{>0})}^n \ni (x_1, \dots, x_n) \freccia u(x_1, \dots, x_n, a_1, \dots, a_m) \in \erre_{>0}$$

For every $t$, this function defines, by conjugation with the map $\log_{\left(\frac{1}{t}\right)}$, a function on $\erre^n$ corresponding to the term $u$ where the operations are interpreted with the operations of $\erre^t$, and every constant $a_i$ is interpreted as $\log_{\left(\frac{1}{t}\right)}(a_i)$:
$$ U_t = \log_{\left(\frac{1}{t}\right)} \circ U \circ {\left(\Log_{\left(\frac{1}{t}\right)}\right)}^{-1}: \erre^n \freccia \erre $$

\begin{lemma}
Let $U_0: \erre^n \freccia \erre$ be the function defined by the term $u$ where the operations are interpreted with the operations of $\erre^{trop}$, and every constant $a_i$ is interpreted as $0$. Then 
$$\forall x \in \erre^n :  U_0(x) \leq U_t(x) \leq U_0(x) + \log_{\left( \frac{1}{t} \right)} C $$
where $C$ is a constant depending only on the term $u$ and the coefficients $a_i$.
In particular the family of functions $U_t$ uniformly converges to the function $U_0$.
\end{lemma}

\begin{proof}
By induction on the complexity of the term. If $u=x_1$, then $U_0 = U_t$ and $C = 1$. If $u = y_1$ then $U_t = \log_{\left( \frac{1}{t} \right)} a_1$ and $U_0 = 0$, hence $C = a_1$. If $u = v^\alpha$, where $\alpha \in \erre$, then $V_0 \leq V_t \leq V_0 + \log_{\left( \frac{1}{t} \right)} C$, hence $U_0 \leq U_t \leq U_0 + \log_{\left( \frac{1}{t} \right)} C^\alpha $. If $u = v \cdot w$, then $V_0 \leq V_t \leq V_0 + \log_{\left( \frac{1}{t} \right)} C$  and  $W_0 \leq W_t \leq W_0 + \log_{\left( \frac{1}{t} \right)} D$, hence $U_0 \leq U_t \leq U_0 + \log_{\left( \frac{1}{t} \right)} CD $. If $u = v + w$, then $V_0 \leq V_t \leq V_0 + \log_{\left( \frac{1}{t} \right)} C$ and $W_0 \leq W_t \leq W_0 + \log_{\left( \frac{1}{t} \right)} D$, hence $U_0 \leq U_t \leq U_0 + \log_{\left( \frac{1}{t} \right)} 2\max(C,D) $.
\end{proof}

If $\phi(x_1, \dots, x_n, y_1, \dots, y_m)$ is an $L_{\mathcal{OS}^\erre}$-formula and $a_1, \dots, a_m$ are constants, they define the set:
$$V = \{ (x_1, \dots, x_n) \in {(\erre_{>0})}^n \ | \ \phi(x_1, \dots, x_n, a_1, \dots, a_m) \}$$

We will denote by $\phi_t$ the formula $\phi$ where the operations are interpreted in the structure $\erre^t$, and $\phi_0$ the formula $\phi$ where the operations are interpreted in the structure $\erre^{trop}$. Hence

$$\ameba_t(V) = \{(x_1, \dots, x_n) \in \erre^n \ | \ \phi_t(x_1, \dots, x_n, \log_{\left( \frac{1}{t} \right)} a_1, \dots, \log_{\left( \frac{1}{t} \right)} a_m) \} $$
Because $\log_{\left(\frac{1}{t}\right)}$ is a semifield isomorphism hence the amoeba $\ameba_t(V)$ is described by the same formula. Anyway it is not always true that
$$\ameba_0(V) = \{(x_1, \dots, x_n) \in \erre^n \ | \ \phi_t(x_1, \dots, x_n, 0, \dots, 0 \}$$
For example if $\phi(x_1) = \neg (x \leq 1)$, then $\phi_0 = \neg ( x \leq 0)$, but the logarithmic limit set of $\{x > 1\}$ is not $\{ x > 0\}$, but $\{ x \geq 0\}$.

\subsection{Dequantization of sets}

A \nuovo{positive formula} is a formula written without the symbols $\neg, \Rightarrow, \Leftrightarrow$. These formulae contains only the connectives $\vee$ and $\wedge$ and the quantifiers $\forall$, $\exists$. Consider the standard $\mathcal{OS}^\erre$-structure on $\erre_{>0}$, or one of the $\mathcal{OS}^\erre$-structures $\erre^t$ or $\erre^{trop}$ on $\erre$. Every subset of ${(\erre_{>0})}^n$ or $\erre^n$ that is defined by a quantifier-free positive $L_{\mathcal{OS}^\erre}$-formula in one of these structures is closed, as the set of symbols $\mathcal{OS}^\erre$ has only the relations $=$ and $\leq$, that are closed, and the functions $+$, $\cdot$, $x^\alpha$, that are continuous.

\begin{proposition}
Let $\phi(x_1, \dots, x_n, y_1, \dots, y_m)$ be a positive $L_{\mathcal{OS}^\erre}$-formula, and let $a_1, \dots, a_m \in \erre_{>0}$ be parameters. If $V$ is such that
$$V \subset \{ (x_1, \dots, x_n) \in {(\erre_{>0})}^n \ |\ \phi(x_1, \dots, x_n, a_1, \dots, a_m) \}$$
Then
$$\ameba_0(V) \subset \{ x \in \erre^n \ |\  \phi_0(x_1, \dots, x_n,0, \dots, 0)  \} $$
\end{proposition}

\begin{proof}
By induction on the complexity of the formula. If $\phi$ is atomic, then it has the form
$$u(x_1, \dots, x_n, y_1, \dots, y_m) \ \mathcal{R}\  v(x_1, \dots, x_n, y_1, \dots, y_m)$$
where $\mathcal{R}$ is $=$ or $\leq$. We have
$$\ameba_t(V) \subset \{ x \in \erre^n \ |\  \phi_t(x_1, \dots, x_n, \log_{\frac{1}{t}}(a_1), \dots, \log_{\frac{1}{t}}(a_m) \}  $$

We may put all the equations together, one for every $t$, thus finding a description for the deformation
$$\mathcal{D} = \{(x,t) \in \erre^n \times (0,\varepsilon) \ |\ x \in \ameba_t(V) \}$$
$$\mathcal{D} = \{(x_1, \dots, x_n,t) \in \erre^n \times (0,1) \ |\ \phi_t(x_1, \dots, x_n, \log_{\frac{1}{t}}(a_1), \dots, \log_{\frac{1}{t}}(a_m) \}$$

If we consider $U_t$ and $V_t$ as functions on $\erre^n \times (0,1)$, they can be extended continuously to $\erre^n \times [0,1)$ defining the extensions on $\erre^n \times \{0\}$ by $U_0, V_0$. Hence we get following inclusion for the logarithmic limit set:
$$\ameba_0(V) \subset \{ x \in \erre^n \ |\  \phi_0(x_1, \dots, x_n,0, \dots, 0)  \} $$

If $\phi = \psi_1 \vee \psi_2$ (resp. $\phi = \psi_1 \wedge \psi_2$), then $V \subset V_1 \cup V_2$ (resp. $V \subset V_1 \cap V_2$), where $V_i$ is defined by $\psi_i$. The statement follows from proposition \ref{prop:loglimsetproperties}.

If $\phi = \exists x_{n+1} : \psi$, then $V$ is contained in the projection of $W$, where $W$ is the set defined by $\psi$. The statement follows from proposition \ref{prop:projection}.

If $\phi = \forall x_{x+1}: \psi$, then we denote by $W$ the set defined by $\psi$. If $(0,\dots,0,-1) \in \ameba_0(V)$, there is a sequence $x(k)$ in $V$ converging to a point $(b_1, \dots, b_{n-1},0)$ with $b_i \neq 0$. Then $W$ contains a sequence of lines $\{(x(k),y) \ |\ y \in \erre_{>0}\}$, hence $\ameba_0(W)$ contains the line $\{(0,\dots,0,y) \ |\ y \in \erre\}$. As $\ameba_0(W) \subset \{(x_1, \dots, x_{n+1}) \ |\ \psi_0(x_1, \dots, x_{n+1})\}$, then $\ameba_0(V) \subset \{(x_1, \dots, x_n) \ |\ \forall x_{n+1} : \psi_0(x_1, \dots, x_{n+1})\}$.
\end{proof}

Anyway there are examples where 
$$V = \{x \in {(\erre_{>0})}^n \ |\ \phi(x_1, \dots, x_n, a_1, \dots, a_m)\}$$ 
with $\phi$ a positive $\mathcal{OS}^\erre$-formula, and
$$\ameba_0(V) \subsetneq  \{x \in \erre^n \ |\ \phi_0(x_1, \dots, x_n, 0, \dots, 0)\}$$

\begin{figure}[htbp]
\begin{center}
\includegraphics[height=\figuresize]{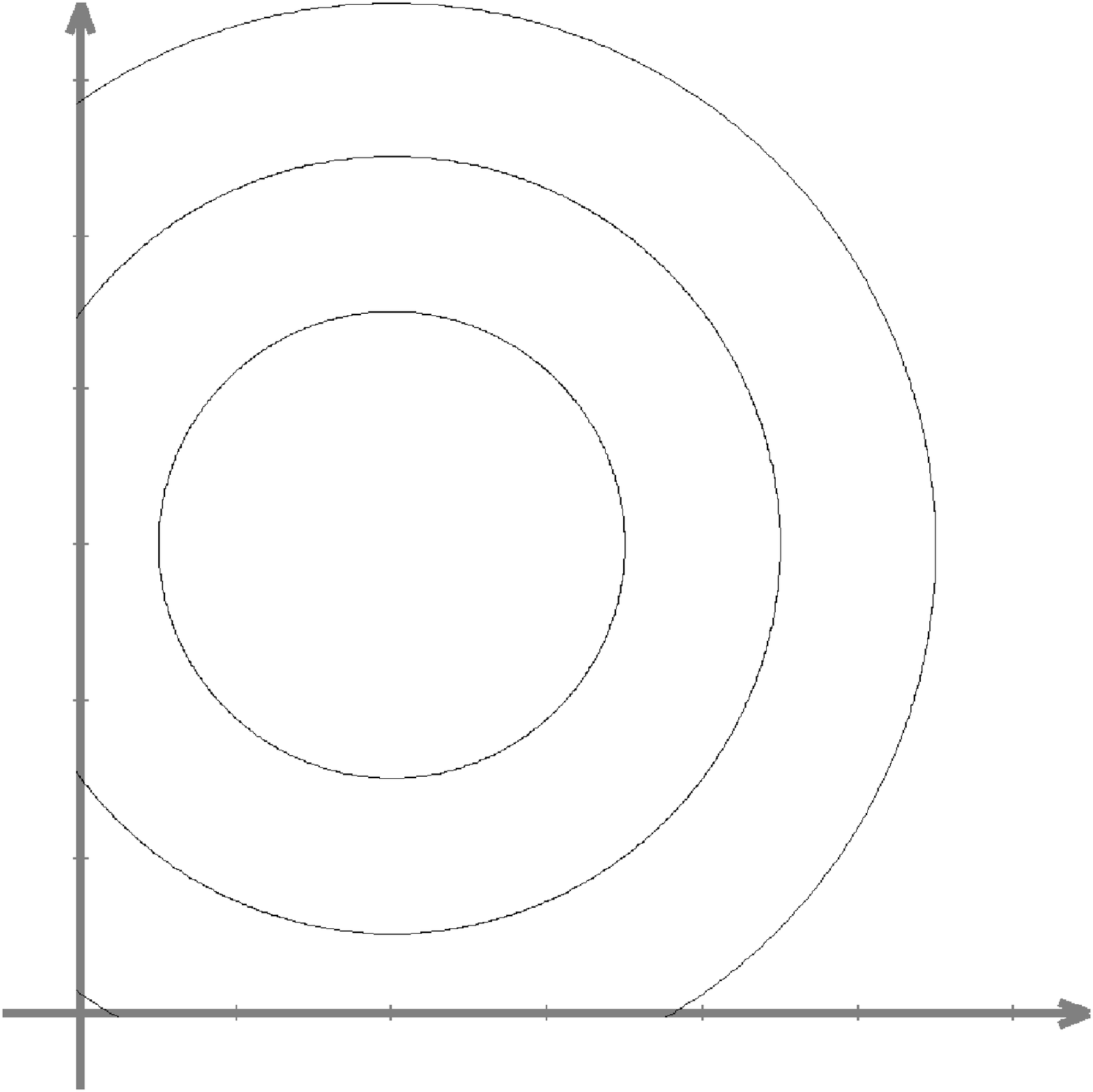}
\caption{$V_r = \{(x,y) \in {(\erre_{>0})}^2 \ |\  x^2 + y^2 + 13 = r^2 + 4 x + 6 y \}$ for $r = \frac{3}{2}, \frac{5}{2}, \frac{7}{2}$ respectively. Their logarithmic limit sets are different.}\label{fig:circle}
\end{center}  
\end{figure}

For example consider the following atomic formula
$$\phi(x_1,x_2,y_1,y_2,y_3): x_1^2 + x_2^2 + y_1 = y_2 x_1 + y_3 x_2$$
with constants $a_1 = 13 - r^2, a_2 = 4, a_3 = 6$, with $r^2 < 13$. This is the equation of a circle,
${(x_1-2)}^2 + {(x_2-3)}^2 = r^2$. The dequantized formula does not depend on the value of $r$:
$$\phi_0(x_1,x_2,0,0,0,0): \max(2x_1,2x_2,0) = \max(x_1,x_2)$$
Now if 
$$V_r = \{(x_1, x_2)\in {(\erre_{>0})}^n \ |\ x_1^2 + x_2^2 + (13 - r^2) = 4 x_1 + 6 x_2\}$$
the logarithmic limit sets of $V_{\frac{3}{2}}$, $V_{\frac{5}{2}}$, $V_{\frac{7}{2}}$ are different (see figure \ref{fig:circle}). We have
$$\ameba_0(V_{\frac{7}{2}}) = \{(x_1, x_2) \in \erre^2 \ |\ \max(2x_1,2x_2,0) = \max(0,x_1,x_2)\}$$
but for $\ameba_0(V_{\frac{3}{2}})$ and $\ameba_0(V_{\frac{5}{2}})$ we have a strict inclusion.

Even if $\phi(x_1, x_2, \tfrac{27}{4}, 4, 6)$ is a definition of $V_{\frac{5}{2}}$, $\phi_0(x_1, x_2,0,0,0)$ is not a definition of $\ameba_0(V_{\frac{5}{2}})$. Anyway we can find another formula with this property, for example:
$$(4x_1^3 + 16x_1^2 + 37x_1 + 40)(x_1^2 -4x_1 + x_2^2 -6x_2 + \frac{27}{4}) = 0$$
The dequantized version of this formula is
$$\max(5x_1,0,2x_2,3x_1+2x_2) = x_2+\max(3x_1, 0) $$
And this formula is an exact description of $\ameba_0(V_{\frac{5}{2}})$.

There are examples of real algebraic sets $V$ where it is not possible to find an algebraic formula $\phi$ defining $V$ such that $\phi_0$ is a definition of $\ameba_0(V)$. Here by an algebraic formula we mean an atomic formula with the relation $=$, in other words an equation between two positive polynomials.

Consider for example the cubic
$$V = \{(x_1,x_2) \in \erre^2 \ |\ x_1^2 + x_2^2 + 1 = 2x_2 + x_1^3\} $$
as in figure \ref{fig:cubic}. This cubic has an isolated point in $(0,1)$. This point is outside the positive orthant ${(\erre_{>0})}^2$, hence it does not influence the logarithmic limit set of $V' = V \cap {(\erre_{>0})}^2$, but the set defined by
$$ \{ (x_1, x_2) \in \erre^2 \ |\ \max(2x_1, 2x_2, 0) = \max(x_2, 3x_1) \} $$
contains also the half line $\{x_2 = 0, x_1 \leq 0\}$ that is not in the logarithmic limit set, and the same happens for every polynomial equation defining $V$.

We need to use the order relation $\leq$ to construct a formula $\phi$ defining $V'$ such that $\phi_0$ is a definition of $\ameba_0(V')$. For example:
$$V' = \{ (x_1,x_2) \in {(\erre_{>0})}^2 \ |\ x_1^2 + x_2^2 + 1 = 2x_2 + x_1^3 \ \wedge\ x_1 \geq \frac{1}{2} \} $$
$$\ameba_0(V') = \{ (x_1,x_2) \in \erre^2 \ |\ \max(2x_1, 2x_2, 0) = \max(x_2, 3x_1) \ \wedge\ x_1 \geq 0 \} $$

As we will see in the next subsection, this is a general fact.

\subsection{Exact definition}

Let $C$ be an open convex set such that $(0,\dots,0,-1) \in C$, and the closure $\overline{C}$ is a convex polyhedral cone contained in $\{x \in \erre^n \ |\ x_n < 0\} \cup \{0\}$. The faces $F_1, \dots, F_k$ of $\overline{C}$ are described by equations
$$a^1_i x_1 + \dots + a^{n-1}_i x_{n-1} + x_n = 0$$
and $C$ is described by
$$C = \{ x \in \erre^n \ |\ x_n < 0 \mbox{ and } \forall i \in \{1, \dots, k\} : a^1_i x_1 + \dots + a^{n-1}_i x_{n-1} + x_n < 0\}$$

For every $h \in \erre_{>0}$, consider the set
$$E_h(C) = \{ x \in {(\erre_{>0})}^n \ |\ x_n < h \mbox{ and } \forall i \in \{1, \dots, k\} : x_1^{a^1_i } \dots  x_{n-1}^{a^{n-1}_i} x_n < h\} $$

\begin{lemma}
Let $V \subset {(\erre_{>0})}^n$ be a set such that $\ameba_0(V) \cap C = \emptyset$. Then for every sufficiently small $h \in \erre_{>0}$ we have $V \cap E_h(C) = \emptyset$.
\end{lemma}

\begin{proof}
Suppose that for all $i \in \enne$ there exists $x_i \in V \cap E_{\frac{1}{i}}(C)$. Then from the sequence $(x_i) \subset V$ we can extract a subsequence $y_i$ such that $\Log_e(y_i)$ converges to a point $y \in C$.
\end{proof}

Note that $E_h(C)$ is described by the following $S_{OS}^\erre$-formula, with $y=h$
$$\phi^C(x_1, \dots, x_n,y) = \neg(y \leq x_n \vee y \leq x_1^{a^1_1 } \dots x_{n-1}^{a^{n-1}_1} x_n   \vee \dots \vee y \leq x_1^{a^1_n } \dots x_{n-1}^{a^{n-1}_n} x_n)$$
and $C$ is described by the formula $\phi^C_0$ with $y=0$.

Let $C \subset \erre^n$ be an open convex set such that the closure $\overline{C}$ is a convex polyhedral cone and $\overline{C} \subset H \cup \{0\}$ where $H$ is an open half-space $H$.

There exists a linear map $B$ such that $(0,\dots,0,-1) \in B(C)$, and $B(\overline{C})$ is contained in $\{x \in \erre^n \ |\ x_n < 0\} \cup \{0\}$. We will use the notation 
$$E_h(C) = \overline{B}^{-1}( E_h(B(C)) )$$
As before there exists a $S_{OS}^\erre$-formula $\phi^C(x_1, \dots, x_n, y)$ such that
$$E_h(C) = \{x \in {(\erre_{>0})}^n \ |\ \phi^C(x_1, \dots, x_n,h)\} $$
$$C = \{x \in \erre^n \ |\ \phi^C_0(x_1, \dots, x_n,0) \}$$

Let $V \subset {(\erre_{>0})}^n$ be a set definable in an o-minimal, polynomially bounded structure with field of exponents $\erre$. Then, by theorem \ref{teo:polyh complex}, $\ameba_0(V)$ is a polyhedral complex, hence we can find a finite number of sets $C_1, \dots, C_k$ such that sets such that 

\begin{enumerate}
 \item $C_1 \cup \dots \cup C_k$ is the complement of $\ameba_0(V)$.
 \item The closure $\overline{C_i}$ is a convex polyhedral cone.
 \item There exists an open half-space $H_i$ such that $\overline{C_i} \subset H_i \cup \{0\}$.
\end{enumerate}

\begin{lemma}    \label{lemma:trop descr}
Consider the $S_{OS}^\erre$-formula
$$\phi(x_1, \dots, x_n,y) = \neg\phi^{C_1}(x_1, \dots, x_n,y) \wedge \dots \wedge \neg\phi^{C_k}(x_1, \dots, x_n,y))$$
Then 
$$\ameba_0(V) = \{x \in \erre^n \ |\ \phi_0(x_1, \dots, x_n,0)\}$$
and for every sufficiently small $h \in \erre_{>0}^n$ we have
$$V \subset \{ x \in {(\erre_{>0})}^n \ |\ \phi(x_1, \dots, x_n,h)\}$$
\end{lemma}

\begin{proof}
The first assertion is trivial, and the second assertion follows from previous lemma.
\end{proof}

Note that the formula $\phi$ of previous lemma has the form $\psi_1 \wedge \dots \wedge \psi_k $,
where the $\psi_i$ have the form:
$$\psi_i = y \leq x_1^{a^1_1 } \dots x_{n}^{a^{n}_1}  \vee \dots \vee y \leq x_1^{a^1_m } \dots x_{n}^{a^{n}_m} $$

These formulae does not contain the $+$ operation, hence when they are interpreted with the dequantizing operations $\oplus_t, \odot_t$ or the tropical operations $\oplus, \odot$ the interpretation does not depend on $t$, and it is simply:
$$\psi_i = y \leq {a^1_1 }x_1 +\dots+ {a^{n}_1}x_{n} \vee \dots \vee y \leq {a^1_m}x_1 + \dots + {a^{n}_m} x_{n} $$

\begin{corollary}
Let $V$ be definable in an o-minimal, polynomially bounded structure with field of exponents $\erre$. For $\varepsilon > 0$, and for small enough $t > 0$:
$$\sup_{x \in \ameba_t(V)} d(x,\ameba_0(V)) < \varepsilon$$
\end{corollary}

\begin{proof}
Choose $h$ such that $V \subset \{\phi(x_1, \dots, x_n,h)\}$. Then $\ameba_t(V) \subset \{\phi_t(x_1, \dots, x_n,\log_{\left(\frac{1}{t}\right)}h)\}$. Note that $\{\phi_t(x_1, \dots, x_n,\log_{\left(\frac{1}{t}\right)}h)\}$ is a uniformly bounded neighborhood of $\ameba_0(V)$, with distance depending linearly on $y$, hence the distance tends to zero when $y$ tends to zero.
\end{proof}

\begin{theorem}
Let $\phi(x_1, \dots, x_n, y_1, \dots, y_m)$ be a positive $L_{\mathcal{OS}^\erre}$-formula, let $a_1, \dots, a_m \in \erre_{>0}$ be parameters and denote
$$V = \{ (x_1, \dots, x_n) \in {(\erre_{>0})}^n \ |\ \phi(x_1, \dots, x_n, a_1, \dots, a_m) \}$$
Then there exists a positive $L_{\mathcal{OS}^\erre}$-formula $\psi(x_1, \dots, x_n, y_1, \dots, y_l)$ and parameters $b_1, \dots, b_l \in \erre_{>0}$ such that:
$$V = \{ (x_1, \dots, x_n) \in {(\erre_{>0})}^n \ |\ \psi(x_1, \dots, x_n, b_1, \dots, b_l) \}$$
$$\ameba_0(V) = \{ x \in \erre^n \ |\  \psi_0(x_1, \dots, x_n,0, \dots, 0)  \} $$
\end{theorem}

\begin{proof}
Let $\phi'(x_1, \dots, x_n,y)$ and $h$ as in lemma \ref{lemma:trop descr}. Then $\psi = \phi \wedge \phi'$ is the searched formula.
\end{proof}

\begin{corollary}
Let $V \subset {(\erre_{>0})}^n$ be a closed semi-algebraic set. Then there exists a positive quantifier-free $L_{\mathcal{OS}^\erre}$-formula $\phi(x_1, \dots, x_n, y_1, \dots, y_m)$ and constants $a_1, \dots, a_m \in \erre_{>0}$ such that
$$V = \{ x \in {(\erre_{>0})}^n \ |\ \phi(x_1, \dots, x_n, a_1, \dots, a_m)\}$$
$$\ameba_0(V) = \{x \in \erre^n \ |\ \phi_0(x_1, \dots, x_n, 0, \dots, 0)\}$$
\end{corollary}

\begin{proof}
By \cite[thm. 2.7.2]{BCR98}, every closed semi-algebraic set is defined by a positive quantifier-free $L_{\mathcal{OS}}$-formula.
\end{proof}

\end{document}